\documentclass[matbio]{svjour}
\usepackage{amsmath,amssymb,times}

\usepackage{bbm,epic}

\newcommand\B{B}
\newcommand\C{\mathrm{C}}
\newcommand\ccint[1]{\mathopen[ #1 \mathclose]}
\newcommand\coint[1]{\mathopen[ #1 \mathclose[}
\newcommand\ooint[1]{\mathopen] #1 \mathclose[}
\newcommand\dist{\operatorname{dist}}
\newcommand\dd[1]{\,\mathrm{d}#1}
\newcommand\eps\varepsilon
\newcommand\Id{\mathbbm{1}}
\newcommand\Ll[1]{\ell^{#1}}
\newcommand\M{\mathcal{M}}
\newcommand\NN{\mathbb{N}}
\newcommand\NNzer{\NN_0}
\newcommand\ord{\operatorname{\scriptstyle\mathcal{O}}}
\newcommand\R{\mathcal{R}}
\newcommand\Ra{\tilde{\R}}
\newcommand\RR{\mathbb{R}}
\newcommand\RRnn{\RR_{\ge0}}
\newcommand\X{X}
\renewcommand\le\leqslant
\renewcommand\ge\geqslant

\allowdisplaybreaks
\newcommand{\es}{\hspace*{0.5pt}}  % space before punctuation marks in numbered equations
\newcommand\is\,  % space between indices, e.g. M_{i \is i+1}
\newcommand\vc[1]{\boldsymbol{#1}}  % vector
\newcommand\vp{\vphantom{()}}

\newcommand\pr[1]{\textup{(}#1\textup{)}}
\newcommand\textdef[1]{\emph{#1}}

\newcommand\dref[1]{Definition~\textup{\ref{def:#1}}}

\newcommand\eref[1]{\textup{(\ref{eq:#1})}}
\newcommand\fref[1]{Figure~\textup{\ref{fig:#1}}}

\newcommand\lref[1]{Lemma~\textup{\ref{lem:#1}}}

\newcommand\pref[1]{Proposition~\textup{\ref{prop:#1}}}

\newcommand\sref[1]{Section~\textup{\ref{sec:#1}}}

\newcommand\thref[1]{Theorem~\textup{\ref{thm:#1}}}

\makeatletter
\newcommand\cf{cf.\@ifnextchar,{}{\ }}
\newcommand\eg{e.g.\@ifnextchar,{}{, }}
\newcommand\ie{i.e.\@ifnextchar,{}{, }}
\makeatother
\newcommand\chaps{Chs.~}
\newcommand\cor{Cor.~}
\newcommand\Eqs{Equations }
\newcommand\lem{Lemma~}
\newcommand\propositions{Propositions }
\newcommand\prop{Prop.~}
\newcommand\remk{Rem.~}
\newcommand\sect{Sec.~}
\newcommand\sections{Sections }
\newcommand\thm{Thm.~}
\newcommand\thms{Thms.~}
\newcommand\weaks{\mbox{weak-$*$} }
\newcommand\Wrt{With respect to }
\newcommand\wrt{with respect to }

\let\oldcite\cite
\def\cite{\upshape\oldcite}
\smartqed
\makeatletter
\newif\if@showqed
\@showqedtrue
\global\@namedef{endproof}{\if@showqed\qed\fi\global\@showqedtrue\@endtheorem}
\makeatother
\makeatletter
\newlength\qedraise
\newcommand\qedhere{\@ifnextchar[{\qed@here}{\qed@here[0pt]}%] bracket matching
}
\def\qed@here[#1]{%
  \global\setlength{\qedraise}{#1}%
  {\@xp\aftergroup\csname\@currenvir @qed\endcsname}\global\@showqedfalse}%
  \def\displaymath@qed{%
    \relax
    \ifmmode
      \ifinner \aftergroup\linebox@qed
      \else
        \eqno
        \let\eqno\relax \let\leqno\relax \let\veqno\relax
        \raisebox{\qedraise}{\qed}%
      \fi
    \else
       \aftergroup\linebox@qed
    \fi
  }
  \@xp\let\csname equation*@qed\endcsname\displaymath@qed
  \@xp\let\csname multline*@qed\endcsname\displaymath@qed
  \def\align@qed{\tag*{\raisebox{\qedraise}{\qed}}}
  \@xp\let\csname align*@qed\endcsname\align@qed
\def\linebox@qed{\hfil\hbox{\qed}\hfilneg}
\makeatother

\begin{document}

\title{Unequal Crossover Dynamics \\ in Discrete and Continuous Time}
\author{Oliver Redner \and Michael Baake}
\institute{Institut f\"ur Mathematik, Univ.\ Greifswald,
  Jahnstr.\ 15a, 17487 Greifswald, Germany 
  \and 
  Fakult\"at f\"ur Mathematik. Univ.\ Bielefeld, Pf.\ 100131, 33501 Bielefeld, Germany
  \and
  \email{redner@uni-greifswald.de}
  \and
  \email{mbaake@mathematik.uni-bielefeld.de}}
\makeatletter
\def\@date{}
\makeatother
\def\copyleft{}
\def\makeheadbox{}
\maketitle

\begin{abstract}
  We analyze a class of models for unequal crossover (UC) of sequences
  containing sections with repeated units that may differ in length.
  In these, the probability of an `imperfect' alignment, in which the
  shorter sequence has $d$ units without a partner in the longer one,
  scales like $q^d$ as compared to `perfect' alignments where all
  these copies are paired.  The class is parameterized by this penalty
  factor $q$.  An effectively infinite population size and thus
  deterministic dynamics is assumed.  For the extreme cases $q=0$ and
  $q=1$, and any initial distribution whose moments satisfy certain
  conditions, we prove the convergence to one of the known fixed
  points, uniquely determined by the mean copy number, in both
  discrete and continuous time.  For the intermediate parameter
  values, the existence of fixed points is shown.
\end{abstract}

\section{Introduction}
\label{sec:intro}

Recombination is a by-product of (sexual) reproduction that leads
to the mixing of parental genes by exchanging genes (or sequence
parts) between homologous chromosomes (or DNA strands). This is
achieved through an alignment of the corresponding sequences,
along with crossover events which lead to reciprocal exchange
of the induced segments. In this process, imperfect alignment
may result in sequences that differ in length form the parental
ones; this is known as {\em unequal crossover\/} (UC). Imperfect
alignment is facilitated by the presence of repeated elements
(as is observed in some rDNA sequences, compare \cite{GrLi}), and 
is believed to be an important driving mechanism for their evolution.
The repeated elements may follow an evolutionary course independent of 
each other and thus give rise to evolutionary innovation. For a detailed
discussion of these topics, see \cite{ShAt,ShWa} and references therein.

This article is concerned with a class of models for UC, originally
investigated by Shpak and Atteson \cite{ShAt} for discrete time, which
is built on preceding work by Ohta \cite{Oht} and Walsh \cite{Wsh}
(see \cite{ShAt} for further references).  Starting from their partly heuristic
results, we prove various existence and uniqueness theorems and
analyze the convergence properties, both in discrete and in continuous
time. This will require a rather careful mathematical development
because the dynamical systems are infinite dimensional.

In this model class, one considers individuals whose genetic
sequences contain a section with repeated units.  These may vary in
number, $i\in\NNzer = \{0,1,2,3,\ldots\}$, where $i=0$ is explicitly
allowed, corresponding to no unit being present (yet).  The
composition of these sections (\wrt mutations that might have
occurred) and the rest of the sequence are ignored here.

In the course of time, recombination events happen, each of which
basically consists of three steps. First, independent pairs are formed at random
(in equidistant time steps, or at a fixed rate). Then, their respective 
sections are randomly aligned, possibly with imperfections in form of `overhangs',
according to some probability distribution for the various possibilities. 
Finally, both sequences are cut at an arbitrary common position between two adjacent
building blocks, with uniform distribution for the cut positions, and
their right (or left) fragments are interchanged.  This so-called
unequal crossover is schematically depicted in \fref{ureco}.
Obviously, the total number of relevant units is conserved in each
event.

\begin{figure}[t]
  \begin{center}
    \unitlength=1mm
     \begin{picture}(85,11)(0,0)
      \matrixput(0.625,9.5)(2.5,0){4}(0,0){1}{\line(1,0){1.25}}
      \matrixput(10,8)(5,0){9}(0,3){2}{\line(1,0){5}}
      \matrixput(10,8)(5,0){10}(0,3){1}{\line(0,1){3}}
      \matrixput(60,8)(5,0){3}(0,3){2}{\line(1,0){5}}
      \matrixput(60,8)(5,0){4}(0,3){1}{\line(0,1){3}}
      \matrixput(75.625,9.5)(2.5,0){4}(0,0){1}{\line(1,0){1.25}}
      \matrixput(0.625,1.5)(2.5,0){10}(0,0){1}{\line(1,0){1.25}}
      \matrixput(25,0)(5,0){6}(0,3){2}{\line(1,0){5}}
      \matrixput(25,0)(5,0){7}(0,3){1}{\line(0,1){3}}
      \matrixput(60,0)(5,0){2}(0,3){2}{\line(1,0){5}}
      \matrixput(60,0)(5,0){3}(0,3){1}{\line(0,1){3}}
      \matrixput(70.625,1.5)(2.5,0){6}(0,0){1}{\line(1,0){1.25}}
      \qbezier(55,9.5)(57.5,5.5)(59.9,1.5)
      \qbezier(55,1.5)(57.5,5.5)(59.9,9.5)
    \end{picture}
    \caption{Snapshot of an unequal crossover event as described in the text.
      Rectangles denote the relevant blocks, while the dashed lines
      indicate possible extensions with other elements that are
      disregarded here.}
    \label{fig:ureco}
  \end{center}
\end{figure}
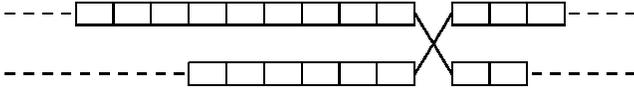

We assume the population size to be (effectively) infinite.  (Concerning
finite populations, see the remarks in \sref{remarks}.)  Then, almost
surely in the probabilistic sense, compare \cite[Sec.\ 11.2]{EtKu86},
the population is described by the deterministic time evolution of
a probability measure $\vc{p} \in \M_1^+(\NNzer)$, which we
identify with an element $\vc{p} = (p_k^{})_{k\in\NNzer}^{}$ in the
appropriate subset of $\Ll{1}(\NNzer)$.  Since we will not consider any
genotype space other than $\NNzer$ in this article, reference to it will be
omitted in what follows.  These spaces are complete in the metric derived from
the usual $\Ll{1}$ norm, which is the same as the total variation norm here.
The metric is denoted by
\begin{equation}
  \label{eq:urgpmetric}
  d(\vc{p},\vc{q}) =
  \|\vc{p} - \vc{q}\|^{}_1 = \sum_{k\ge0} |p_k^{} - q_k^{}| \es.
\end{equation}

With this notation, the above process is described by the recombinator
\begin{equation}
  \label{eq:reco}
  \R(\vc{p})_i^{} = 
  \frac{1}{\|\vc{p}\|^{}_1}
    \sum_{j,k,\ell\ge0} \, T_{ij,k\ell}^{} \, p_k^{} \, p_\ell^{} \es.
\end{equation}
Here, $T_{ij,k\ell} \ge 0$ denotes the probability that, given a pair
$(k,\ell)$, this pair turns into $(i,j)$. Consequently, for normalization, we require
\begin{equation}
  \label{eq:sumt}
  \sum_{i,j\ge0} T_{ij,k\ell} = 1 
  \qquad\text{for all $k, \ell \in \NNzer$.}
\end{equation}
The factor $p_k^{} \, p_\ell^{}$ in \eref{reco} describes the
probability that a pair $(k,\ell)$ is formed, \ie we assume that two
individuals are chosen independently from the population.  We assume
further that $T_{ij,k\ell} = T_{ji,k\ell} = T_{ij,\ell k}$, \ie that
$T_{ij,k\ell}$ is symmetric \wrt both index pairs, which is
reasonable.  Then, the summation over $j$ in \eref{reco}
represents the breaking-up
of the pairs after the recombination event.  These two ingredients of
the dynamics constitute what is known as (\textdef{instant})
\textdef{mixing} and are responsible for the quadratic nature of the
iteration process.

As mentioned above, we will only consider processes that conserve the
total copy number in each event, \ie $T^{(q)}_{ij,k\ell} > 0$ for $i+j
= k+\ell$ only.  To\-geth\-er with the normalization \eref{sumt} and
the symmetry condition from above, this yields the (otherwise weaker) 
condition
\begin{equation}
  \label{eq:urcons}
  \sum_{i\ge0} i \, T_{ij,k\ell}
  = \sum_{i,j\ge0} \frac{i+j}{2} \, T_{ij,k\ell}
  = \frac{k+\ell}{2} \es,
\end{equation}
which implies conservation of the \emph{mean copy number} in the population,
\begin{displaymath}
  \sum_{i\ge0} i \, \R(\vc{p})_i^{}
  = \sum_{i,j,k,\ell\ge0} i \, T^{(q)}_{ij,k\ell} \, p_k^{} \, p_\ell^{}
  = \sum_{k,\ell\ge0} \frac{k+\ell}{2} \, p_k^{} \, p_\ell^{}
  = \sum_{k\ge0} k \, p_k^{} \es.
\end{displaymath}

Condition \eref{sumt} and the presence of the prefactor
$1/\|\vc{p}\|^{}_1$ in \eref{reco} make $\R$ norm non-increasing, \ie
$\|\R(\vc{x})\|^{}_1 \le \|\vc{x}\|^{}_1$, and positive homogeneous of
degree 1, \ie $\R(a\vc{x}) = |a| \R(\vc{x})$, for all $\vc{x} \in \Ll{1}$
and $a \in \RR$.  Furthermore, $\R$ is a positive operator with
$\|\R(\vc{x})\|^{}_1 = \|\vc{x}\|^{}_1$ for all positive elements
$\vc{x} \in \Ll{1}$.  Thus, it is guaranteed that $\R$ maps $\M_r^+$,
the space of positive measures of total mass $r$, into itself.
This space is complete in the
topology induced by the norm $\|.\|_1^{}$, \ie by the metric $d$ from
\eref{urgpmetric}.  (For $r=1$, of course, the prefactor
is redundant but ensures numerical stability of an iteration with
$\R$.)

Given an initial configuration $\vc{p}_0^{} = \vc{p}(t\mathbin=0)$,
the dynamics may be taken in discrete time steps, with subsequent
generations,
\begin{equation}
  \label{eq:urdisctime}
  \vc{p}(t+1) = \R(\vc{p}(t)) \es,
  \qquad t\in\NNzer \es.
\end{equation}
Our treatment of this case will be set up in a way that also allows
for a generalization of the results to the analogous process in
continuous time, where generations are overlapping,
\begin{equation}
  \label{eq:urconttime}
  \tfrac{\dd{}}{\dd{t}} \vc{p}(t) = \varrho \, (\R-\Id)(\vc{p}(t)) \es,
  \qquad t\in\RRnn \es.
\end{equation}
Obviously, the (positive) parameter $\varrho$ in \eref{urconttime} only leads
to a rescaling of the time $t$.  We therefore choose $\varrho=1$ without loss
of generality.  Furthermore, the formulation of discrete versus continuous
time dynamics in \eref{urdisctime} and \eref{urconttime} is chosen
so that the fixed points of \eref{urdisctime} are identical to the equilibria 
of \eref{urconttime}, regardless of $\varrho$. This can easily be verified by
a direct calculation. In what follows, we will thus use the term \emph{fixed 
point} for both discrete and continuous dynamics.

In the UC model, one distinguishes `perfect' alignments, in which each
unit in the shorter sequence has a partner in the longer sequence, and
`imperfect' alignments, with `overhangs' of the shorter sequence relative
to the longer one. To come to a reasonable probability distribution for
the various possiblities, the first are taken to be equally probable among
each other, while the latter are penalized by a
factor $q^d$ relative to the first, where $q \in \ccint{0,1}$ is a
model parameter and $d$ is the length of the overhang (at most the entire
length of the shorter sequence; in the example of \fref{ureco}, interpreted
as a snapshot right after the crossover event took place, we
have $d=1$).  In the extreme case $q=0$, only perfect alignments may
occur, whereas for $q=1$ overhangs are not penalized at all and one
obtains the uniform distribution on the possibilities.  For
obvious reasons, the first case is dubbed \textdef{internal UC}, the
second \textdef{random UC} \cite{ShAt}.

It is now straightforward, though a bit tedious, to derive the transition
probabilities $T^{(q)}_{ij,k\ell}$. To this end, one has to trace what
happens in steps two and three of the recombination event only, while
the random formation of pairs does not enter here.
This has been done in \cite{ShAt} and need not be repeated.
However, in view of our above remarks, it is desirable to rewrite
the findings in a way that reflects the natural symmetry properties
of the $T^{(q)}_{ij,k\ell}$. In compact notation, this leads to the 
transition probabilities
\begin{equation}
  \label{eq:urtq}
  T^{(q)}_{ij,k\ell} = 
  C^{(q)}_{k\ell} \, \delta_{i+j,k+\ell} \, (1+\min\{k,\ell,i,j\}) \,
    q^{0 \vee (k\wedge\ell - i \wedge j)} \es,
\end{equation}
where $k \vee \ell := \max\{k,\ell\}$, $k \wedge \ell :=
\min\{k,\ell\}$, and $0^0 = 1$.  The $C^{(q)}_{k\ell}$ are chosen
such that \eref{sumt} holds, \ie $\sum_{i,j\ge0} T^{(q)}_{ij,k\ell} =
1$, and are hence symmetric in $k$ and $\ell$.  Explicitly, they read
(see also \cite[\sect 2.1]{ShAt})
\begin{displaymath}
  C^{(q)}_{k\ell} =
  \frac{(1-q)^2}{
    (k\wedge\ell+1)(|k-\ell|+1)(1-q)^2 +
    2q(k\wedge\ell - (k\wedge\ell+1)q + q^{k\wedge\ell+1})} \es.
\end{displaymath}
Note further that the total number of units is indeed conserved in
each event and that the process is symmetric within both pairs.  Hence
\eref{urcons} is satisfied.

Let us briefly come back to the question of `discrete' versus `continuous' 
time, which are considered simultaneously for good reasons.
Common to both is the nonlinearity that stems from the probability that a
certain (random) pair is formed in the first place. Then,
for the discrete time dynamics \eref{urdisctime}, the $T^{(q)}_{ij,k\ell}$
have the direct meaning of the transition probability that, given a pair
$(k,\ell)$, this turns into a pair $(i,j)$. In contrast, for the
continuous time dynamics \eqref{eq:urconttime}, the number $T^{(q)}_{ij,k\ell}$
is to be understood as the probability to obtain a pair $(i,j)$
\emph{conditioned} on a recombination event with a pair of type
$(k,\ell)$, of which each recombines at the same rate. In probabilistic
terminology, the $\R$ of \eref{urdisctime} is the discrete time
skeleton of the process in continuous time, also called the
embedded discrete time process.

\bigskip\noindent
The aim of this article is to find answers to the following questions:
\begin{enumerate}
\renewcommand\labelenumi{\theenumi.}
\item Are there fixed points of the dynamics?
\item Given the mean copy number $m$, is there a unique fixed point?
\item If so, under which conditions and in which sense does an initial
  distribution converge to this fixed point under time evolution?
\end{enumerate}
Of course, the trivial fixed point with $p_0^{} = 1$ and $p_k^{} = 0$
for $k>0$ always exists, which we generally exclude from our
considerations.  But even then, the answer to the first question is
positive for general operators of the form \eref{reco} that satisfy
\eref{sumt} and some rather natural further condition.  This is
discussed in \sref{urgeneral}.  For the extreme cases $q=0$ (internal
UC) and $q=1$ (random UC), fixed points are known explicitly for every
$m$ and it has been conjectured \cite{ShAt} that, under mild
conditions, also questions 2 and 3 can be answered positively for all
values of $q \in \ccint{0,1}$.  Indeed, for both extreme cases, norm
convergence of the population distribution to the fixed points can be
shown, which is done in \sections \ref{sec:urq0} and \ref{sec:urq1},
respectively.  Since the dynamical systems involved are infinite
dimensional, a careful analysis of compactness properties is needed
for rigorous answers.  The proofs for $q=1$ are based on alternative
representations of probability measures via generating functions,
presented in \sref{urrepr}.  For the intermediate parameter regime, we
can only show that there exists a fixed point for every $m$, but
neither its uniqueness nor convergence to it, see \sref{urqgen}.  Some
remarks in \sref{remarks} conclude this article.

\section{Existence of fixed points}
\label{sec:urgeneral}

Let us begin by stating the following general fact.
\begin{proposition}
  \label{prop:recolipschitz}
  If the recombinator\/ $\R$ of \eref{reco} satisfies \eref{sumt},
  then the global Lipschitz condition
  \begin{displaymath}
    \|\R(\vc{x}) - \R(\vc{y})\|^{}_1 \le C \|\vc{x} - \vc{y}\|^{}_1
  \end{displaymath}
  is satisfied, with constant\/ $C=3$ on\/ $\Ll{1}$, respectively\/
  $C=2$ if\/ $\vc{x}$, $\vc{y} \in \M_r$.
\end{proposition}
\begin{proof}
  Let $\vc{x}$, $\vc{y} \in \Ll{1}$ be non-zero (otherwise the
  statement is trivial). Then,
  \begin{align*}
    \|\R(\vc{x}) - \R(\vc{y})\|^{}_1 &=
    \sum_{i\ge0} \left|
      \sum_{j,k,\ell\ge0} T_{ij,k\ell}^{} \left(
        \frac{x_k^{}\,x_\ell^{}}{\|\vc{x}\|^{}_1} -
        \frac{y_k^{}\,y_\ell^{}}{\|\vc{y}\|^{}_1} \right) \right| \\
    &\le \sum_{k,\ell\ge0} \left|
      \frac{x_k^{}\,x_\ell^{}}{\|\vc{x}\|^{}_1} -
      \frac{y_k^{}\,y_\ell^{}}{\|\vc{y}\|^{}_1} \right|
    \sum_{i,j\ge0} T_{ij,k\ell} \\
    &= \sum_{k,\ell\ge0} \left|
      \frac{x_k^{}\,x_\ell^{}}{\|\vc{x}\|^{}_1} -
      \frac{x_k^{}\,y_\ell^{}}{\|\vc{x}\|^{}_1} +
      \frac{x_k^{}\,y_\ell^{}}{\|\vc{x}\|^{}_1} -
      \frac{y_k^{}\,y_\ell^{}}{\|\vc{y}\|^{}_1} \right| \\
    &\le \sum_{k,\ell\ge0} \left(
      \frac{|x_k^{}|}{\|\vc{x}\|^{}_1} |x_\ell^{} - y_\ell^{}| +
      |y_\ell^{}| \left| \frac{x_k^{}}{\|\vc{x}\|^{}_1} - 
        \frac{y_k^{}}{\|\vc{y}\|^{}_1} \right| \right) \\
    &= \|\vc{x}-\vc{y}\|^{}_1 + 
    \frac{1}{\|\vc{x}\|^{}_1} \bigl\|
      \|\vc{y}\|^{}_1 \vc{x} - \|\vc{x}\|^{}_1 \vc{y} \bigr\|^{}_1 \es.
  \end{align*}
  The last term becomes
  \begin{multline*}
    \frac{1}{\|\vc{x}\|^{}_1} \bigl\|
      \|\vc{y}\|^{}_1 \vc{x} - \|\vc{x}\|^{}_1 \vc{y} \bigr\|^{}_1 \\
    = \frac{1}{\|\vc{x}\|^{}_1} \bigl\| 
      (\|\vc{y}\|^{}_1 - \|\vc{x}\|^{}_1) \vc{x} +
      \|\vc{x}\|^{}_1 (\vc{x} - \vc{y}) \bigr\|^{}_1
    \le 2 \|\vc{x}-\vc{y}\|^{}_1 \es,
  \end{multline*}
  from which $\|\R(\vc{x}) - \R(\vc{y})\|^{}_1 \le 3\|\vc{x} -
  \vc{y}\|^{}_1$ follows for $\vc{x}$, $\vc{y} \in \Ll{1}$.  If
  $\vc{x}$, $\vc{y} \in \M_r$, the above calculation simplifies to
  $\|\R(\vc{x}) - \R(\vc{y})\|^{}_1 \le 2\|\vc{x} - \vc{y}\|^{}_1$.
\end{proof}
In continuous time, this is a sufficient condition for the existence
and uniqueness of a solution of the initial value problem
\eref{urconttime}, \cf \cite[\thms 7.6 and 10.3]{Ama}.  Another useful
notion in this respect is the following.
\begin{definition}{\cite[\sect 18]{Ama}}
  \label{def:lyapunov}
  Let\/ $Y$ be an open subset of a Banach space\/ $E$ and let\/ $f
  \colon Y \to E$ satisfy a \pr{local} Lipschitz condition.  A
  continuous function\/ $L$ from\/ $X \subset Y$ to\/ $\RR$ is called
  a \textdef{Lyapunov function} for the initial value problem
  \begin{displaymath}
    \tfrac{\dd{}}{\dd{t}} \vc{x}(t) = f(\vc{x}(t)) \es,
    \qquad \vc{x}(0) = \vc{x}_0 \in X \es,
  \end{displaymath}
  if the \textdef{orbital derivative}\/ $\dot{L}(\vc{x}_0) :=
  \liminf_{t\to 0^+} \frac1t \bigl( L(\vc{x}(t)) - L(\vc{x}_0) \bigr)$
  satisfies
  \begin{equation}
    \label{eq:dotl}
    \dot{L}(\vc{x}_0) \le 0
  \end{equation}
  for all initial conditions\/ $\vc{x}_0 \in X$.
\end{definition}
If further $\dot{L}(\vc{x}_{\mathrm{F}}) = 0$ for a single fixed point
$\vc{x}_{\mathrm{F}}$ only, then $L$ is called a \textdef{strict}
Lyapunov function.  If a Lyapunov function exists, we have
\begin{theorem}{\cite[\thm 17.2 and \cor 18.4]{Ama}}
  \label{thm:lyapunov}
  With the notation of \dref{lyapunov}, assume that there is a
  Lyapunov function\/ $L$, that the set\/ $X$ is closed, and that, for
  an initial condition\/ $\vc{x}_0 \in X$, the set\/ $\{ \vc{x}(t) :
  t\in\RRnn, \text{$\vc{x}(t)$ exists} \}$ is relatively compact in\/
  $X$.  Then,\/ $\vc{x}(t)$ exists for all\/ $t\ge0$ and
  \begin{displaymath}
    \lim_{t\to\infty} \dist(\vc{x}(t), X_L) = 0 \es,
  \end{displaymath}
  where\/ $\dist(\vc{x},X_L) = \inf_{\vc{y} \in X_L} \|\vc{x} -
  \vc{y}\|$ and\/ $X_L$ denotes the largest invariant subset of\/ $\{
  \vc{x} \in X : \dot{L}(\vc{x}) = 0 \}$ \pr{in forward and backward
    time}.  \hspace*{\fill}\qed
\end{theorem}
Obviously, if $L$ is a strict Lyapunov function, we have $X_L =
\{\vc{x}_{\mathrm{F}}\}$ and this theorem implies $d(\vc{x}(t),
\vc{x}_{\mathrm{F}}) \to 0$ as $t \to \infty$.

Returning to the original question of the existence of fixed points,
we now recall the following facts, compare \cite{BilC,Shir} for
details.
\begin{proposition}{\cite[\cor to \thm V.1.5]{Yos}}
  \label{prop:vaguenorm}
  Assume the sequence\/ $\bigl(\vc{p}^{(n)}\bigr)$ in\/ $\M_1^+$ to
  converge in the \weaks topology \pr{\ie pointwise, or vaguely} to
  some\/ $\vc{p} \in \M_1^+$, \ie
  \begin{displaymath}
    \lim_{n\to\infty} p^{(n)}_k = p_k^{}
    \quad\text{for all\, $k\in\NNzer$} \es,
    \quad\text{with\, $p_k^{}\ge0$ \; and\; 
      $\textstyle\sum_{k\ge0} p_k^{} = 1$} \es.
  \end{displaymath}
  Then, it also converges weakly \pr{in the probabilistic sense} and in
  total variation, \ie $\lim_{n\to\infty} \|\vc{p}^{(n)} -
  \vc{p}\|^{}_1 = 0$.\hspace*{\fill}\qed
\end{proposition}
\begin{proposition}
  \label{prop:urgfixed}
  Assume that the recombinator\/ $\R$ from \eref{reco} satisfies
  \eref{sumt} and has a convex, \weaks closed invariant set\/ $M
  \subset \M_1^+$, \ie $\R(M) \subset M$, that is\/ \emph{tight}, \ie for
  every\/ $\eps>0$ there is an\/ $m\in\NNzer$ such that\/ $\sum_{k \ge
    m} p_k^{} < \eps$ for every\/ $\vc{p} \in M$.  Then, $\R$ has a
  fixed point in\/ $M$.
\end{proposition}
\begin{proof}
  Prohorov's theorem \cite[\thm III.2.1]{Shir} states that tightness
  and relative compactness in the \weaks topology are equivalent (see
  also \cite[\chaps 1.1 and 1.5]{BilC}).  In our case, $M$ is tight
  and \weaks closed, therefore, due to \pref{vaguenorm}, norm compact.
  Furthermore, $M$ is convex and $\R$ is (norm) continuous by
  \pref{recolipschitz}.  Thus, the claim follows from the
  Leray--Schauder--Tychonov fixed point theorem \cite[\thm
  V.19]{ReSi}.
\end{proof}
\Wrt the UC model, we will see that such compact invariant subsets
indeed exist.

\section{Internal unequal crossover}
\label{sec:urq0}

After these preliminaries, let us begin with the case of internal UC
with perfect alignment only, \ie $q=0$ in \eref{urtq}.  This case is
the simplest because, in each recombination event, no sequences exceeding
the longer of the participating sequences can be formed. Here, on $\M_1^+$,
the recombinator \eref{reco} simplifies to
\begin{equation}
  \label{eq:ur0reco}
  \R_0(\vc{p})_i^{} =
  \sum_{\substack{k,\ell\ge0 \\[0.4\baselineskip] 
      k\wedge\ell \le i \le k\vee\ell}} 
    \frac{p_k^{} \, p_\ell^{}}{1+|k-\ell|} \es.
\end{equation}
From now on, we write $\R_q$ rather than $\R$ whenever we look at a
recombinator with (fixed) parameter $q$.  It is instructive to
generalize the notion of reversibility (or detailed balance, compare
\cite[(4.1)]{ShAt}).
\begin{definition}
  We call a probability measure\/ $\vc{p} \in \M_1^+$
  \textdef{reversible} for a recombinator\/ $\R$ of the form
  \eref{reco} if, for all\/ $i, j, k, \ell \ge 0$,
  \begin{equation}
    \label{eq:urrev}
    T_{ij,k\ell}^{} \, p_k^{} \, p_\ell^{} = 
    T_{k\ell,ij}^{} \, p_i^{} \, p_j^{} \es.
  \end{equation}
\end{definition}
The relevance of this concept is evident from the following property.
\begin{lemma}
  \label{lem:urrev}
  If\/ $\vc{p} \in \M_1^+$ is reversible for\/ $\R$, it is also a
  fixed point of\/ $\R$.
\end{lemma}
\begin{proof}
  Assume $\vc{p}$ to be reversible.  Then, by \eref{sumt},
  \begin{displaymath}
    \R(\vc{p})_i^{} =
    \sum_{j,k,\ell\ge0} \, T_{ij,k\ell}^{} \, p_k^{} \, p_\ell^{} =
    \sum_{j,k,\ell\ge0} \, T_{k\ell,ij}^{} \, p_i^{} \, p_j^{} = 
    p_i^{} \sum_{j\ge0} p_j^{} = p_i^{} \es. \qedhere[-13.9pt]
  \end{displaymath}
\end{proof}
So, in our search for fixed points, we start by looking for solutions
of \eref{urrev}.  Since, for $q=0$, forward and backward transition
probabilities are simultaneously non-zero only if $\{i,j\} =
\{k,\ell\} \subset \{n,n+1\}$ for some $n$, the components $p_k^{}$
may only be positive on this small set as well.  By the following
proposition, this indeed characterizes all fixed points for $q=0$.
\begin{proposition}
  \label{prop:ur0unique}
  A probability measure\/ $\vc{p} \in \M_1^+$ is a fixed point of\/
  $\R_0$ if and only if its mean copy number\/ $m = \sum_{k\ge0}
  k\,p_k^{}$ is finite,\/ $p_{\lfloor m \rfloor}^{} = \lfloor m
  \rfloor + 1 - m$, \linebreak $p_{\lceil m \rceil}^{} = m + 1 -
  \lceil m \rceil$, and\/ $p_k^{} = 0$ for all other\/ $k$.  This
  includes the case that\/ $m$ is integer and\/ $p_{\lfloor m
    \rfloor}^{} = p_{\lceil m \rceil}^{} = p_m^{} = 1$.
\end{proposition}
\begin{proof}
  The `if' part was stated in \cite[\sect 4.1]{ShAt} and follows
  easily by insertion into \eref{ur0reco} or \eref{urrev}.  For the
  `only if' part, let $i$ denote the smallest integer such that
  $p_i^{} > 0$.  Then,
  \begin{displaymath}
    \R(\vc{p})_i^{} = 
    p_i^2 + 2 p_i^{} \sum_{\ell\ge1} \frac{p_{i+\ell}^{}}{1+\ell} =
    p_i^{} \left( p_i^{} + p_{i+1}^{} + 
      \sum_{\ell\ge2} \frac{2}{\ell+1} p_{i+\ell}^{} \right) \le 
    p_i^{} \es,
  \end{displaymath}
  where the last step follows since $\frac{2}{\ell+1} < 1$ in the last
  sum, with equality if and only if $p_k^{} = 0$ for all $k \ge i+2$.
  This implies $m < \infty$ and the uniqueness of $\vc{p}$ (given $m$)
  with the non-zero frequencies as claimed.
\end{proof}

It it possible to analyze the case of internal UC on the basis of the
compact sets to be introduced below in \sref{urrepr}.  However, as J.
Hofbauer pointed out to us \cite{HofPriv}, it is more natural to start
with a larger compact set to be introduced in \eref{m1mc}.  Our main
result in this section is thus
\begin{theorem}
  \label{thm:ur0conv}
  Assume that, for the initial condition\/ $\vc{p}(0)$ and fixed\/
  $r>1$, the\/ $r$-th moment exists,\/ $\sum_{k\ge0} k^r p_k^{}(0) <
  \infty$.  Then, $m = \sum_{k\ge0} k \, p_k^{}(0)$ is finite and,
  both in discrete and in continuous time,\/ $\lim_{t\to\infty}
  \|\vc{p}(t) - \vc{p}\|^{}_1 = 0$ with the appropriate fixed point\/
  $\vc{p}$ from \pref{ur0unique}.
\end{theorem}
The proof relies on the following lemma, which slightly modifies and
completes the convergence arguments of \cite[\sect 4.1]{ShAt}, puts
them on rigorous grounds, and extends them to continuous time.
\begin{lemma}
  \label{lem:ur0dp}
  Let\/ $r>1$ be arbitrary, but fixed.  Consider the set of
  probability measures with fixed mean\/ $m<\infty$ and a centered\/
  $r$-th moment bounded by\/ $C<\infty\es$,
  \begin{equation}
    \label{eq:m1mc}
    \M_{1,m,C}^+ = \{ \vc{p}\in\M_1^+ : \sum_{k\ge0} k\,p_k^{} = m,\,
      M_r(\vc{p}) \le C \} \es,
  \end{equation}
  equipped with \pr{the metric induced by} the total variation norm,
  where
  \begin{equation}
    \label{eq:ur0dp}
    M_s(\vc{p}) = \sum_{k\ge0} |k-m|^s \, p_k^{}
  \end{equation}
  for\/ $s \in \{1,r\}$.  This is a compact and convex space.  Both\/
  $M_1$ and\/ $M_r$ satisfy\/ $M_s(\R_0(\vc{p})) \le M_s(\vc{p})$,
  with equality if and only if\/ $\vc{p}$ is a fixed point of\/
  $\R_0$.  Furthermore,\/ $M_1$ is a continuous mapping from\/
  $\M_{1,m,C}^+$ to\/ $\RRnn$ and a Lyapunov function for the dynamics
  in continuous time.
\end{lemma}
\begin{proof}
  Let a sequence $(\vc{p}^{(n)}) \subset \M_{1,m,C}^+$ be given and
  consider the random variables $\vc{f}^{(n)} = (k)_{k\in\NNzer}^{}$
  on the probability spaces $(\NNzer,\vc{p}^{(n)})$.  Their
  expectation values are equal to $m$, which, by Markov's inequality
  \cite[p.~599]{Shir}, implies the tightness of the sequence
  $(\vc{p}^{(n)})$.  Hence, by Prohorov's theorem \cite[\thm
  III.2.1]{Shir} (see also \cite[\chaps 1.1 and 1.5]{BilC}), it
  contains a convergent subsequence $(\vc{p}^{(n_i)})$ (recall that,
  by \pref{vaguenorm}, norm and pointwise convergence are equivalent
  in this case).  Let $\widetilde{\vc{p}} \in \M_1^+$ denote its limit and
  $\widetilde{\vc{f}} = (k)_{k\in\NNzer}^{}$ a random variable on
  $(\NNzer,\widetilde{\vc{p}})$, to which the $\vc{f}^{(n_i)}$ converge in
  distribution.  Since $r>1$, the $\vc{f}^{(n_i)}$ are uniformly
  integrable by Markov's and H\"older's inequalities.  Hence, due to
  \cite[\lem 3.11]{Kal}, their expectation values, which all equal
  $m$, converge to the one of $\widetilde{\vc{f}}$, which is thus $m$ as
  well.  Consider now the random variables $\vc{g}^{(n_i)} =
  \widetilde{\vc{g}} = (|k-m|^r)_{k\in\NNzer}^{}$ on
  $(\NNzer,\vc{p}^{(n)})$ and $(\NNzer,\widetilde{\vc{p}})$, respectively.
  The expectation values of the $\vc{g}^{(n_i)}$ are bounded by $C$,
  which, again by \cite[\lem 3.11]{Kal}, is then also an upper bound
  for the expectation value of $\widetilde{\vc{g}}$ (to which the
  $\vc{g}^{(n_i)}$ converge in distribution).  This proves the
  compactness of $\M_{1,m,C}^+$.  The convexity is obvious.

  \Wrt the second statement, consider
  \begin{equation}
    \label{eq:ur0dr0}
    \begin{aligned}
      M_s(\R_0(\vc{p})) &=
      \sum_{i\ge0} \!
        \sum_{\substack{k,\ell\ge0 \\[0.4\baselineskip] 
            k\wedge\ell \le i \le k\vee\ell}} \!
        \frac{|i-m|^s}{1+|k-\ell|} \, p_k^{} \, p_\ell^{} \\
      &= \sum_{k,\ell\ge0} \frac{p_k^{} \, p_\ell^{}}{1+|k-\ell|}
        \, \frac12 \sum_{i=k\wedge\ell}^{k\vee\ell}
        (|i-m|^s + |k+\ell-i-m|^s) \es.
    \end{aligned}
  \end{equation}
  For notational convenience, let $j = k+\ell-i$.  We now show
  \begin{equation}
    \label{eq:ur0ijklm}
    |i-m|^s + |k+\ell-i-m|^s \le |k-m|^s + |\ell-m|^s \es.
  \end{equation}
  If $\{k,\ell\} = \{i,j\}$, then \eref{ur0ijklm} holds with equality.
  Otherwise, assume, without loss of generality, that $k < i \le j <
  \ell$.  If $m \le k$ or $m \ge \ell$, we have equality for $s=1$ but
  a strict inequality for $s=r$ due to the convexity of $x \mapsto
  x^r$.  (For $s=1$, this describes the fact that a recombination
  event between two sequences that are both longer or both shorter
  than the mean does not change their mean distance to the mean copy
  number.)  In the remaining cases, the inequality is strict as well.
  Hence, $M_s(\R_0(\vc{p})) \le M_s(\vc{p})$ with equality if and only
  if $\vc{p}$ is a fixed point of $\R_0$, since otherwise the sum in
  \eref{ur0dr0} contains at least one term for which \eref{ur0ijklm}
  holds as a strict inequality.
  
  To see that $M_1$ is continuous, select a converging sequence
  $(\vc{p}^{(n)})$ in $\M_{1,m,C}^+$ and consider the random variables
  $\vc{h}^{(n)} = (|k-m|)_{k\in\NNzer}^{}$ on $(\NNzer,\vc{p}^{(n)})$.
  As above, the latter are uniformly integrable, from which the
  continuity of $M_1$ follows.  Since $M_1(\vc{p})$ is linear in
  $\vc{p}$ and thus infinitely differentiable, so is the solution
  $\vc{p}(t)$ for every initial condition $\vc{p}_0 \in \M_{1,m,C}^+$,
  compare \cite[\thm 9.5 and \remk 9.6(b)]{Ama}.  Therefore, we have
  \begin{displaymath}
    \dot{M_1}(\vc{p}_0) =
    \liminf_{t\to0^+} \frac{M_1(\vc{p}(t)) - M_1(\vc{p}_0)}{t} =
    M_1(\R_0(\vc{p}_0)) - M_1(\vc{p}_0) \le 0 \es,
  \end{displaymath}
  again with equality if and only if $\vc{p}_0$ is a fixed point.
  Thus, $M_1$ is a Lyapunov function.
\end{proof}
\begin{theopargself}
  \begin{proof}[of \thref{ur0conv}]
    By assumption, the $r$-th moment of $\vc{p}(0)$ exists, which is
    equivalent to the existence of the centered $r$-th moment by
    Minkowski's inequality \cite[\sect II.6.6]{Shir}.  This obviously
    implies the existence of the mean $m$.  By \lref{ur0dp},
    $\vc{p}(t) \in \M_{1,m,C}^+$ follows for all $t\ge0$,
    directly for discrete time and via a satisfied subtangent
    condition \cite[\thm VI.2.1]{Mar} (see also \cite[\thm 16.5]{Ama})
    for continuous time.  In the discrete case, due to the compactness
    of $\M_{1,m,C}^+$, there is a convergent subsequence
    $(\vc{p}(t_i))$ with some limit $\vc{p}$.  Consider now the mean
    distance $M_1$ to the mean copy number from \eref{ur0dp}.  If
    $\lim_{t\to\infty} \vc{p}(t) = \vc{p}$, we have, due to the
    continuity of $M_1$ and $\R_0$,
    \begin{displaymath}
      M_1(\R_0(\vc{p})) = \lim_{t\to\infty} M_1(\R_0(\vc{p}(t))) =
      \lim_{t\to\infty} M_1(\vc{p}(t+1)) = M_1(\vc{p}) \es,
    \end{displaymath}
    thus $\vc{p}$ is a fixed point by \lref{ur0dp}.  Otherwise, there
    are two convergent subsequences $(\vc{p}(t_i))$, with limit
    $\vc{p}$, and $(\vc{p}(s_i))$, with limit $\vc{q}$, whose indices
    alternate, $t_i < s_i < t_{i+1}$.  Then, we also have
    $M_1(\R_0(\vc{p}(t_i))) \ge M_1(\vc{p}(s_i))$ and
    $M_1(\R_0(\vc{p}(s_i))) \ge M_1(\vc{p}(t_{i+1}))$, and therefore
    \begin{multline*}
      M_1(\vc{p}) \ge M_1(\R_0(\vc{p})) = 
      \lim_{i\to\infty} M_1(\R_0(\vc{p}(t_i))) \ge
      \lim_{i\to\infty} M_1(\vc{p}(s_i))
      = M_1(\vc{q}) \\
      \ge M_1(\R_0(\vc{q})) = 
      \lim_{i\to\infty} M_1(\R_0(\vc{p}(s_i))) \ge
      \lim_{i\to\infty} M_1(\vc{p}(t_{i+1})) = M_1(\vc{p}) \es.
    \end{multline*}
    Thus, both $\vc{p}$ and $\vc{q}$ are fixed points by \lref{ur0dp}
    and hence equal by \pref{ur0unique}.  In continuous time, the
    claim follows from \thref{lyapunov} since $M_1$ is a Lyapunov
    function by \lref{ur0dp}.
  \end{proof}
\end{theopargself}

Note that, for $q=0$, the recombinator can be expressed in terms of
explicit frequencies $\pi_{k,\ell}$ of fragment pairs before
concatenation (with copy numbers $k$ and $\ell$) as $\R_0(\vc{p})_i^{}
= \sum_{j=0}^i \pi_{j,i-j}$.  However, we have, so far, not been able
to use this for a simplification of the above treatment.

\section{Alternative probability representations}
\label{sec:urrepr}

Our next goal is to find the analogue of \thref{ur0conv} for the case of
$q=1$ (random UC).
Whereas the convergence arguments for the case $q=0$ relied on a compact set
of probability measure defined via the $r$-th moment, we are not (yet) able to
extend this approach to $q>0$.  Instead, we will consider, as an alternative
representation for a probability measure $\vc{p} \in \M_1^+$, the generating
function
\begin{equation}
  \label{eq:urgpsip}
  \psi(z) = \sum_{k\ge0} p_k^{} z^k \es,
\end{equation}
for which $\psi(1) = \|\vc{p}\|^{}_1 = 1$ and the radius of
convergence is at least 1.  We will restrict our discussion to such
$\vc{p}$ for which $\limsup_{k\to\infty} \sqrt[k]{p_k^{}} < 1$.  Then,
the radius of convergence, $\rho(\psi) = 1/\limsup_{k\to\infty}
\sqrt[k]{p_k^{}}$ by Hadamard's formula \cite[10.5]{RudRCA}, is larger
than $1$.  This is, biologically, no restriction since for any
`realistic' system there are only finitely many non-zero $p_k^{}$ (and
therefore $\rho(\psi) = \infty$).  Mathematically, this condition
ensures the existence of all moments and enables us to go back and
forth between the probability measure and its generating function,
even when looked at $\psi(z)$ only in the vicinity of $z=1$ (see
\pref{urgainxad} below and \cite[\sect II.12]{Shir}).  By abuse of
notation, we define the induced recombinator for these generating
functions as
\begin{equation}
  \label{eq:urgrecopsi}
  \R(\psi)(z) = \sum_{k\ge0} \R(\vc{p})_k^{} \, z^k \es.
\end{equation}
In general, with the exception of the case $q=1$, we do not know any
simple expression for $\R(\psi)$ in terms of $\psi$.  Nevertheless,
\eref{urgrecopsi} will be central to our further analysis.

It is advantageous to use the local expansion around $z=1$, written in
the form
\begin{equation}
  \label{eq:urgpsi1}
  \psi(z) = \sum_{k\ge0} (k+1) a_k (z-1)^k \es,
\end{equation}
whose coefficients are given by
\begin{equation}
  \label{eq:urga}
  a_k = \frac{1}{(k+1)!} \psi^{(k)}(1) =
  \frac{1}{k+1} \sum_{\ell \ge k} \binom{\ell}{k} \, p_\ell^{} =: 
  \vc{a}(\vc{p})_k^{} \ge 0 \es.
\end{equation}
In particular, $a_0 = 1$ and $a_1 = \frac12 \sum_{\ell\ge0} \ell \,
p_\ell^{}$.  This definition of $a_k$ is size biased, and will become
clear from the simplified dynamics for $q=1$ that results from it.
For the sake of compact notation, we use $\vc{a} =
(a_k)^{}_{k\in\NNzer}$ both for the coefficients and for the mapping.
The coefficients $\vc{a}$ are elements of the following compact,
convex metric space.
\begin{definition}
  \label{def:xad}
  For fixed\/ $\alpha$ and\/ $\delta$ with\/ $0 < \alpha \le \delta <
  \infty$, let
  \begin{displaymath}
    \X_{\alpha,\delta} = 
    \{ \vc{a} = (a_k)_{k\in\NNzer} : a_0 = 1,\, a_1 = \alpha,\,
      \text{$0 \le a_k \le \delta^k$ for $k\ge2$} \} \es.
  \end{displaymath}
  On this space, define the metric
  \begin{equation}
    \label{eq:urgmetric}
    d(\vc{a},\vc{b}) = \sum_{k\ge0} d_k |a_k-b_k|
  \end{equation}
  with\/ $d_k = (\gamma/\delta)^k$ for some\/ $0<\gamma<\frac13$.
\end{definition}
It is obvious that $d$ is indeed a metric and that
$\X_{\alpha,\delta}$ is a convex set, \ie we have $\eta\,\vc{a} +
(1-\eta)\vc{b} \in \X_{\alpha,\delta}$ for all $\vc{a}$, $\vc{b} \in
\X_{\alpha,\delta}$ and $\eta \in \ccint{0,1}$.  Note that we use the
same symbol $d$ as in \eref{urgpmetric} since it will always be clear
which metric is meant.  The space $\X_{\alpha,\delta}$ is naturally
embedded in the Banach space (\cf \cite[\sect 24.I]{WalEn})
\begin{equation}
  \label{eq:hdelta}
  H_{\gamma/\delta} = \{ \vc{x}\in\RR^{\NNzer} : \|\vc{x}\| < \infty \}
\end{equation}
with the norm $\|\vc{x}\| = \sum_{k\ge0} (\gamma/\delta)^k |x_k|$, for
$\gamma$ and $\delta$ as in \dref{xad}.  In particular,
$d(\vc{a},\vc{b}) = \|\vc{a}-\vc{b}\|$.  Furthermore, we have the
following two propositions.
\begin{proposition}
  \label{prop:urgxadcomp}
  The space\/ $\X_{\alpha,\delta}$ is compact in the metric\/ $d$ of
  \eref{urgmetric}.
\end{proposition}
\begin{proof}
  In metric spaces, compactness and sequential compactness are
  equivalent, compare \cite[\thm II.3.8]{Lan}.  Hence, let
  $(\vc{a}^{(n)})$ be any sequence in $\X_{\alpha,\delta}$.  By
  assumption, $a^{(n)}_0 \equiv 1$ and $a^{(n)}_1 \equiv \alpha$.
  Furthermore, each element sequence $(a_k^{(n)}) \subset
  \ccint{0,\delta^k}$ has a convergent subsequence.  We now
  inductively define, for every $k$, a convergent subsequence
  $(a_k^{(n_{k,i})})$, with limit $a_k$, such that the indices
  $\{n_{k,i} : i\in\NN\}$ are a subset of the preceding indices
  $\{n_{k-1,i} : i\in\NN\}$.  This way, we can proceed to a `diagonal'
  sequence $(\vc{a}^{(n_{i,i})})$.  The latter is now shown to
  converge to $\vc{a} = (a_k)$, which is obviously an element of
  $\X_{\alpha,\delta}$.  To this end, let $\eps>0$ be given.  Choose
  $m$ large enough such that $\sum_{k>m} (2\gamma)^k < \eps/2$, and
  then $i$ such that $\sum_{k=0}^m d_k^{} |a^{(n_{i,i})}_k - a_k^{}| <
  \eps/2$.  Then
  \begin{equation}
    \label{eq:urgxadcomp}
    \begin{aligned}
      d(\vc{a}^{(n_{i,i})},\vc{a}) &= 
      \sum_{k=2}^m d_k |a_k^{(n_{i,i})}-a_k^{}| + 
        \sum_{k>m} d_k |a_k^{(n_{i,i})}-a_k^{}| \\
      &< \frac\eps2 + \sum_{k>m} (2\gamma)^k <
      \eps \es,
    \end{aligned}
  \end{equation}
  which proves the claim.
\end{proof}
\begin{proposition}
  \label{prop:urgainxad}
  If\/ $\limsup_{k\to\infty} \sqrt[k]{p_k^{}} < 1$, the coefficients\/
  $a_k$ from \eref{urga} exist and\/ $\vc{a}(\vc{p}) \in
  \X_{\alpha,\delta}$ with\/ $\alpha = \vc{a}(\vc{p})^{}_1 = \frac12 m
  = \frac12 \sum_{k\ge0} k\,p_k^{}$ and some\/ $\delta$.  Conversely,
  if\/ $\vc{p}(\vc{a}) \in \X_{\alpha,\delta}$ for some\/ $\alpha$,
  $\delta$, one has\/ $\limsup_{k\to\infty} \sqrt[k]{p_k^{}} < 1$.
\end{proposition}
For a proof, we need the following
\begin{lemma}
  \label{lem:urgradii}
  Let\/ $f_0(z) = \sum_{k\ge0} c_k z^k$ be a power series with
  non-negative coefficients\/ $c_k$ and\/ $f_x(z) = \sum_{k\ge0}
  \frac{1}{k!} f_0^{(k)}(x) (z-x)^k$ the expansion of\/ $f_0$ around
  some\/ $x \in \coint{0,\rho(f_0)}$.  Then, $\rho(f_0) = x +
  \rho(f_x)$, including the case that both radii of convergence are
  infinite.
\end{lemma}
\begin{proof}
  Since the open disc $\B_x(\rho(f_0)-x)$ is entirely included in
  $\B_0(\rho(f_0))$, the inequality $\rho(f_x) \ge \rho(f_0)-x$
  immediately follows from the theorem of representability by power
  series \cite[\thm 10.16]{RudRCA}.  Consider now the power series
  $f_{xe^{i\varphi}}(z) = \sum_{k\ge0} \frac{1}{k!}
  f_0^{(k)}(xe^{i\varphi}) (z-xe^{i\varphi})^k$ with arbitrary
  $\varphi\in\coint{0,2\pi}$.  Its coefficients satisfy
  $|f_0^{(k)}(xe^{i\varphi})| \le \sum_{n \ge k} \frac{n!}{(n-k)!}
  c_k x^{n-k} = f_0^{(k)}(x)$ due to the non-negativity of the $c_k$.
  This implies $\rho(f_{xe^{i\varphi}}) \ge \rho(f_x)$ by Hadamard's
  formula.  Therefore, $f_0$ admits an analytic continuation on
  $\B_0(x+\rho(f_x))$, the uniqueness of which follows from the
  monodromy theorem \cite[\thm 16.16]{RudRCA}.  The theorem of
  representability by power series then implies the inequality
  $\rho(f_0) \ge x+\rho(f_x)$, which, together with the opposite
  inequality above, proves the claim.
\end{proof}
\begin{theopargself}
  \begin{proof}[of \pref{urgainxad}]
    The assumption implies $\rho(\psi) > 1$ for $\psi$ from
    \eref{urgpsip}.  Then, from \lref{urgradii}, we know that
    $\limsup_{k\to\infty} \sqrt[k]{(k+1)a_k} < \infty$.  Since
    furthermore $a_k \le (k+1) a_k$, also $\limsup_{k\to\infty}
    \sqrt[k]{a_k} < \infty$, so there is an upper bound $\delta$ for
    $\sqrt[k]{a_k}$ and thus $\vc{a}(\vc{p}) \in \X_{\alpha,\delta}$.
    The converse statement follows from \eref{urga} and
    \lref{urgradii}.
  \end{proof}
\end{theopargself}
Therefore, any mapping from $\X_{\alpha,\delta}$ into itself that is
continuous \wrt the metric $d$ from \eref{urgmetric} has a fixed point
by the Leray--Schauder--Tychonov theorem \cite[\thm V.19]{ReSi}.

Note further that each $\X_{\alpha,\delta}$ contains a maximal element
\wrt the partial order $\vc{a}\le\vc{b}$ defined by $a_k \le b_k$ for
all $k\in\NNzer$, which is given by
$(1,\alpha,\delta^2,\delta^3,\ldots)$.  This property finally leads to
\begin{proposition}
  \label{prop:urgpadcomp}
  The space\/ $P_{\alpha,\delta} := \{ \vc{p}\in\M_1^+ :
  \vc{a}(\vc{p})\in\X_{\alpha,\delta} \}$, equipped with \pr{the
    metric induced by} the total variation norm, is compact and
  convex.
\end{proposition}
The proof is based on the following two lemmas.
\begin{lemma}
  \label{lem:urglocbound}
  For any subset of\/ $P_{\alpha,\delta}$, the corresponding
  generating functions from \eref{urgpsip} are locally bounded on\/
  $\B_{1+1/\delta}(0)$.
\end{lemma}
\begin{proof}
  It is sufficient to show boundedness on every compact $K \subset
  \B_{1+1/\delta}(0)$, see \cite[\sect 7.1]{Rem}.  Thus, let such a
  $K$ be given and fix $r \in \coint{0,\frac1\delta}$ so that $K$ is
  contained in $\overline{\B_{1+r}(0)}$.  Then, for every $\vc{p} \in
  P_{\alpha,\delta}$ and every $z \in K$,
  \begin{align*}
    |\psi(z)| &= \Bigl|\sum_{k\ge0} p_k^{} z^k\Bigr| \le
    \sum_{k\ge0} p_k^{} (1+r)^k = \psi(1+r) \\
    &= \sum_{k\ge0} (k+1) \vc{a}(\vc{p})_k^{} \, r^k \le
    1 + 2\,\alpha\,r + \sum_{k\ge2} (k+1) (r\delta)^k <
    \infty \es,
  \end{align*}
  where $r\delta < 1$ was used. This needed to be shown.
\end{proof}
\begin{lemma}
  \label{lem:urgpnormconv}
  If, for a sequence\/ $(\vc{p}^{(n)}) \subset P_{\alpha,\delta}$, the
  coefficients\/ $\vc{a}^{(n)} = \vc{a}(\vc{p}^{(n)})$ from
  \eref{urga} converge to some\/ $\vc{a}$ \wrt the metric\/ $d$ from
  \eref{urgmetric}, then the generating functions\/ $\psi_n$ from
  \eref{urgpsip} converge compactly to some\/ $\psi$ with\/ $\psi(z) =
  \sum_{k\ge0} p_k^{} z^k$ and the\/ $\vc{p}^{(n)}$ thus converge in
  norm to\/ $\vc{p} \in P_{\alpha,\delta}$.
\end{lemma}
\begin{proof}
  By to \lref{urglocbound}, the sequence $(\psi_n)$ is locally
  bounded in $\B_{1+1/\delta}(0)$.  Due to the pointwise convergence
  $|a^{(n)}_k - a^{\vp}_k| \le d_k^{-1} d(\vc{a}^{(n)},
  \vc{a}) \to 0$, we have
  \begin{displaymath}
    \psi^{(k)}_n(1) = (k+1)! \, a^{(n)}_k 
    \;\xrightarrow{n\to\infty}\;
    (k+1)! \, a^{\vp}_k = \psi^{(k)}(1)
  \end{displaymath}
  for every $k\in\NNzer$.  Then, the compact convergence $\psi_n \to
  \psi$ follows from Vitali's theorem \cite[\thm 7.3.2]{Rem}.  In
  particular, this implies that $p^{(n)}_k \to p_k^{\vp} \ge 0$ and $1
  = \sum_{k\ge0} p^{(n)}_k = \psi_n(1) \to \psi(1) = \sum_{k\ge0}
  p_k^{\vp}$, thus $\vc{p} \in \M_1^+$.
  
  Now, choose $r \in \ooint{1,1+\frac1\delta}$.  Then there is, for
  every $\eps>0$, an $n_\eps$ such that $\sup_{|z| \le r}
  |\psi(z)-\psi_n(z)| < \eps$ for all $n \ge n_\eps$.  This implies
  \begin{displaymath}
    |p^{(n)}_k - p_k^{\vp}| =
    \left|\frac{1}{2\pi i} \oint_{|z|=r} 
      \frac{\psi_n(z) - \psi(z)}{z^{k+1}} \dd{z} \right| < \frac{\eps}{r^k}
  \end{displaymath}
  for all $n \ge n_\eps$ by Cauchy's integral formula \cite[\thm
  7.3]{LanCA}.  Now, let $\eps>0$ be given.  Then
  \begin{displaymath}
    \|\vc{p}^{(n)} - \vc{p}\|^{}_1 = 
    \sum_{k\ge0} |p^{(n)}_k - p_k^{\vp}| <
    \eps \frac{1}{1-\tfrac1r}
  \end{displaymath}
  for all $n \ge n_\eps$, which proves the claim.
\end{proof}
\begin{theopargself}
  \begin{proof}[of \pref{urgpadcomp}]
    Let $(\vc{p}^{(n)})$ denote an arbitrary sequence in
    $P_{\alpha,\delta}$ and \linebreak $(\vc{a}^{(n)}) =
    (\vc{a}(\vc{p}^{(n)}))$ the corresponding sequence in
    $\X_{\alpha,\delta}$.  Due to \pref{urgxadcomp}, there is a
    convergent subsequence $(\vc{a}^{(n_i)})^{}_i$.  Then, by
    \lref{urgpnormconv}, $(\vc{p}^{(n_i)})^{}_i$ converges in norm to
    some $\vc{p} \in P_{\alpha,\delta}$.  This proves the compactness
    property.  The convexity of $P_{\alpha,\delta}$ is a simple
    consequence of the convexity of $\M_1^+$, the linearity of the
    mapping $\vc{a}$, and the convexity of $\X_{\alpha,\delta}$.
  \end{proof}
\end{theopargself}

Another property of the mapping $\vc{a} \colon P_{\alpha,\delta} \to
\X_{\alpha,\delta}$ is stated in
\begin{lemma}
  \label{lem:urgacont}
  For every\/ $\alpha$ and\/ $\delta$, the mapping\/ $\vc{a} \colon
  P_{\alpha,\delta} \to \X_{\alpha,\delta}$ from \eref{urga} is
  continuous \pr{\wrt the total variation norm and the metric\/ $d$}
  and injective.  Its inverse\/ $\vc{p} \colon
  \vc{a}(P_{\alpha,\delta}) \to P_{\alpha,\delta}$ is continuous as
  well.
\end{lemma}
\begin{proof}
  Let $\vc{p}, \vc{q} \in P_{\alpha,\delta}$ and assume
  $\vc{a}(\vc{p}) = \vc{a}(\vc{q})$.  Then, as in the proof of
  \lref{urgradii}, the uniqueness of the generating function in
  $\B_{1+1/\delta}(0)$ follows, and thus $\vc{p} = \vc{q}$, which
  proves the injectivity of $\vc{a}$.  The other statements follow
  from Vitali's theorem \cite[\thm 7.3.2]{Rem}:  Norm convergence of a
  sequence $(\vc{p}^{(n)})$ in $P_{\alpha,\delta}$ to some $\vc{p}$
  implies convergence of its element sequences and thus compact
  convergence of the corresponding generating functions $\psi_n$ to
  $\psi$, which is given by $\psi(z) = \sum_{k\ge0} p_k^{}\,z^k$.
  This, in turn, implies convergence of each sequence
  $(\vc{a}(\vc{p}^{(n)})_k)$ to $\vc{a}(\vc{p})_k^{}$, from which, as
  in \eref{urgxadcomp}, the convergence $(\vc{a}(\vc{p}^{(n)})) \to
  \vc{a}(\vc{p})$ (\wrt $d$) follows.  The converse is the statement
  of \lref{urgpnormconv} (see also \cite[\prop 1.6.8]{Ped}).
\end{proof}
Note that, if $\rho(\psi)>2$, the inverse of the mapping $\vc{a}$ is
given by
\begin{displaymath}
  \vc{p}(\vc{a})_k^{}
  = \sum_{\ell \ge k} (-1)^{\ell-k} \binom{\ell}{k} (\ell+1) \, a_\ell^{} \es.
\end{displaymath}

\section{Random unequal crossover}
\label{sec:urq1}

Let us now turn to the random UC model, described by $q=1$ in
\eref{urtq}.  Here, the recombinator \eref{reco} simplifies to
\cite[(3.1)]{ShAt}
\begin{equation}
  \label{eq:ur1reco}
  \R_1(\vc{p})_i^{} =
  \sum_{\substack{k,\ell\ge0 \\[0.4\baselineskip] k+\ell \ge i}}
    \frac{1+\min\{k,\ell,i,k+\ell-i\}}{(k+1)(\ell+1)} 
    \, p_k^{} \, p_\ell^{} \es.
\end{equation}
As for internal UC, by \lref{urrev}, the reversibility condition
\eref{urrev} directly leads to an expression for fixed points,
\begin{displaymath}
  \frac{p_{k\vphantom{j}}^{}}{k+1} \frac{p_{\ell\vphantom{j}}^{}}{\ell+1} =
  \frac{p_{i\vphantom{j}}^{}}{i+1} \frac{p_j^{}}{j+1}
  \quad\text{for all } k+\ell = i+j \es.
\end{displaymath}
This has $p_k^{} = C (k+1) x^k$ as a solution, with appropriate
parameter $x$ and normalization constant $C$.  Again, it turns out
that all fixed points are given this way.
\begin{proposition}{\cite[\thm A.2]{ShAt}}
  \label{prop:ur1fixed}
  Every fixed point\/ $\vc{p} \in \M_1^+$ of\/ $\R_1$ is of the form
  \begin{equation}
    \label{eq:ur1fixed}
    p_k^{} =
    \left(\frac{2}{m+2}\right)^2 (k+1) \left(\frac{m}{m+2}\right)^k \es,
  \end{equation}
  where\/ $m = \sum_{k\ge0} k \, p_k^{} \ge 0$.
  \hspace*{\fill}\qed
\end{proposition}
One can verify this in several ways, one being a direct inductive
calculation.

The main result of this section is
\begin{theorem}
  \label{thm:ur1conv}
  Assume that\/ $\limsup_{k\to\infty} \sqrt[k]{p_k^{}(0)} < 1$.  Then,
  both in discrete and in continuous time,\/ $\lim_{t\to\infty}
  \|\vc{p}(t) - \vc{p}\|^{}_1 = 0$ with the appropriate fixed point\/
  $\vc{p}$ from \pref{ur1fixed}.
\end{theorem}
For a proof, we consider the following alternative process, verbally
described in \cite[p.\ 720f]{ShAt}. It is a two-step stick breaking
and glueing procedure which ultimately induces the same (deterministic) 
dynamics as random UC, even though the underlying process is rather different.  
This will lead to a simple expression for the
induced recombinator of the coefficients $\vc{a}$ from \eref{urga},
which allows for an explicit solution.
\begin{proposition}
  Let\/ $\vc{p} \in \M_1^+$.  Then,
  \begin{equation}
    \label{eq:ur1iterqp}
    \pi_k^{} = \sum_{\ell \ge k} \frac{1}{\ell+1} p_\ell^{}
  \end{equation}
  gives a probability measure\/ $\vc{\pi} \in \M_1^+$, and the
  recombinator can be written as
  \begin{equation}
    \label{eq:ur1iterpq}
    \R_1(\vc{p})_i^{} = 
    \sum_{j=0}^i \pi_j^{} \, \pi_{i-j}^{} = 
    (\vc{\pi}*\vc{\pi})_i^{} \es,
  \end{equation}
  where\/ $*$ denotes the convolution in\/ $\Ll{1}(\NNzer)$.
\end{proposition}
Here, \eref{ur1iterqp} describes a breaking process in which, without any pairing, 
each sequence is cut equally likely between any two of its building blocks. In a
second step, described by \eref{ur1iterpq}, these fragments are paired
randomly and joined (or `glued').
\begin{proof}
  It is easily verified that $\vc{\pi}$ is normalized to $1$.  \Wrt
  \eref{ur1iterpq}, note the following identity for $k+\ell \ge i$,
  \begin{displaymath}
    |\{j : (i-\ell) \vee 0 \le j \le i \wedge k\}| =
    1+\min\{k,\ell,i,k+\ell-i\} \es,
  \end{displaymath}
  which can be shown by treating the four cases on the LHS separately.
  With this, inserting \eref{ur1iterqp} into the RHS of
  \eref{ur1iterpq} yields
  \begin{align*}
    \sum_{j=0}^i \pi_j^{} \, \pi_{i-j}^{} &=
    \sum_{j=0}^i \sum_{k \ge j} \sum_{\ell \ge i-j}
      \frac{1}{(k+1)(\ell+1)} \, p_k^{} \, p_\ell^{} \\
    &= \sum_{\substack{k,\ell \ge 0 \\[0.4\baselineskip] k+\ell \ge i}}
      \frac{1}{(k+1)(\ell+1)} \, p_k^{} \, p_\ell^{}
    \sum_{j=(i-\ell)\vee0}^{i \wedge k} 1 \\
    &= \sum_{\substack{k,\ell\ge0 \\[0.4\baselineskip] k+\ell \ge i}}
      \frac{1+\min\{k,\ell,i,k+\ell-i\}}{(k+1)(\ell+1)} 
      \, p_k^{} \, p_\ell^{} \,=\,
    \R_1(\vc{p})_i^{} \es. \qedhere[-23.4pt]
  \end{align*}
\end{proof}

This nice structure has an analogue on the level of the generating
functions.
\begin{proposition}
  Under the assumptions of \thref{ur1conv},
  let\/ $\phi(z) = \sum_{k\ge0} \pi_k z^k$ denote the generating
  function for\/ $\vc{\pi}$ from \eref{ur1iterqp}. Then,
  \begin{displaymath}
    \phi(z) = \frac{1}{1-z} \int_z^1 \psi(\zeta) \dd{\zeta}
    \qquad\text{and}\qquad
    \R_1(\psi)(z) = \phi(z)^2 \es.
  \end{displaymath}
\end{proposition}
\begin{proof}
  Recall that $\psi(z) = \sum_{k\ge 0} p_k z^k$.
  \Eqs \eref{ur1iterqp} and \eref{ur1iterpq} lead to
  \begin{align*}
    \phi(z) &= 
    \sum_{k\ge0} \sum_{\ell \ge k} \frac{1}{\ell+1} p_\ell^{} z^k =
    \sum_{\ell\ge0} \frac{1}{\ell+1} p_\ell^{} \sum_{k\le\ell} z^k =
    \sum_{\ell\ge0} \frac{1}{\ell+1} p_\ell^{} \frac{1-z^{\ell+1}}{1-z}\\
    &= \frac{1}{1-z} \sum_{\ell\ge0} p_\ell^{} \frac{1-z^{\ell+1}}{\ell+1} =
    \frac{1}{1-z} \sum_{\ell\ge0} p_\ell^{} \int_z^1 \zeta^\ell \dd{\zeta} =
    \frac{1}{1-z} \int_z^1 \psi(\zeta) \dd{\zeta}
  \end{align*}
  and, due to absolute convergence of the series involved,
  \begin{align*}
    \R_1(\psi)(z) &= \sum_{k\ge0} \R_1(\vc{p})_k^{} z^k =
    \sum_{k\ge0} z^k \sum_{\ell=0}^k \pi_\ell \pi_{k-\ell} \\ &=
    \sum_{\ell\ge0} \pi_\ell z^\ell \sum_{k\ge\ell} \pi_{k-\ell} z^{k-\ell} =
    \phi(z)^2 \es.  \qedhere[-13.9pt]
  \end{align*}
\end{proof}
The following lemma states that the radius of convergence of $\psi$
does not decrease under the random UC dynamics.  Thus, it is ensured
that, if $\rho(\psi) > 1$, also $\R_1(\psi)$ may be described by an
expansion at $z=1$, \ie by coefficients $\vc{a}$.
\begin{lemma}
  \label{lem:ur1rhotildepsi}
  The radius of convergence of\/ $\R_1(\psi)$ is\/ $\rho(\R_1(\psi))
  \ge \rho(\psi)$.
\end{lemma}
\begin{proof}
  As $1/\rho(\psi) = \limsup_{k\to\infty} \sqrt[k]{p_k^{}} =: x \le 1$
  and $\lim_{k\to\infty} \sqrt[k]{k+1} = 1$, there is a constant $C>0$
  with $p_k^{} \le C (k+1) x^k$ for all $k$.  Note the identity
  \begin{displaymath}
    \sum_{j=0}^n (1+\min\{i,j,n-i,n-j\}) =(i+1)(n-i+1)
  \end{displaymath}
  for $i \le n$, which follows from an elementary calculation. Then,
  \eref{ur1reco} implies
  \begin{align*}
    \R_1(\vc{p})_i^{} &\le 
    C^2 \sum_{\substack{k,\ell\ge0 \\ k+\ell \ge i}} 
      x^{k+\ell} (1+\min\{k,\ell,i,k+\ell-i\}) \\
    &= C^2 \sum_{n \ge i} x^n \sum_{j=0}^n (1+\min\{i,j,n-i,n-j\}) \\
    &= C^2 (i+1) x^i \sum_{\ell\ge0} (\ell+1) \, x^\ell = 
    \left(\frac{C}{1-x}\right)^2 (i+1) \, x^i \es.
  \end{align*}
  Accordingly, $\limsup_{k\to\infty} \sqrt[k]{\R_1(\vc{p})_k^{}} \le x
  \le 1$, which proves the claim.
\end{proof}

These results enable us to derive the following expression for the
coefficients $\vc{a}$, using the expansion of \eref{urgpsi1}:
\begin{equation}
  \label{eq:ur1recopsi}
  \begin{aligned}
    \R_1(\psi)(z) &= 
    \left[ \frac{1}{z-1} \int_1^z \psi(\zeta) \dd{\zeta} \right]^2 \\
    &= \left[ \sum_{k\ge0} a_k (z-1)^k \right]^2 =
    \sum_{k\ge0} \left( \sum_{n=0}^k a_n a_{k-n} \right) (z-1)^k \es.
  \end{aligned}
\end{equation}
So it is natural to define the induced
recombinator
\begin{equation}
  \label{eq:ur1recoa}
  \Ra_1(\vc{a})_k^{} = 
  \frac{1}{k+1} \sum_{n=0}^k a_n a_{k-n} \ge 0 \es,
\end{equation}
for which we have
\begin{lemma}
  \label{lem:ur1cont}
  The recombinator\/ $\Ra_1$ given by \eref{ur1recoa} maps each
  space\/ $\X_{\alpha,\delta}$ into itself and is continuous \wrt the
  metric\/ $d$ from \eref{urgmetric}.
\end{lemma}
\begin{proof}
  Let $\alpha, \delta > 0$ be given and $\vc{a}, \vc{b} \in
  \X_{\alpha,\delta}$.  Trivially, $\Ra_1(\vc{a})_0^{} = 1$ and
  $\Ra_1(\vc{a})_1^{} = \alpha$.  For $k\ge2$, $\Ra_1(\vc{a})_k^{} =
  \frac{1}{k+1} \sum_{\ell \le k} a_\ell \, a_{k-\ell} \le \delta^k$.
  This proves the first statement.  For the continuity, note first
  that every $\Ra_1(\vc{a})_k^{}$ with $k\ge2$ is continuous as a
  mapping from $\X_{\alpha,\delta}$ to $\ccint{0,\delta^k}$.  Now, let
  $\eps>0$ be given.  Choose $n$ large enough so that $\sum_{k>n}
  (2\gamma)^k < \eps/2$, where $\gamma$ is the parameter introduced in
  \dref{xad}.  Then, there is an $\eta>0$ such that $\sum_{k=2}^n
  (\gamma/\delta)^k |\Ra_1(\vc{a})_k - \Ra_1(\vc{b})_k| < \eps/2$ for
  $\vc{a}, \vc{b} \in \X_{\alpha,\delta}$ with $d(\vc{a},\vc{b}) <
  \eta$.  Thus, for such $\vc{a}$ and $\vc{b}$,
  \begin{displaymath}
    d(\Ra_1(\vc{a}), \Ra_1(\vc{b})) \le
    \sum_{k=0}^n \left(\frac\gamma\delta\right)^k 
      |\Ra_1(\vc{a})_k - \Ra_1(\vc{b})_k| +
    \sum_{k>n} (2\gamma)^k < \eps \es,
  \end{displaymath}
  which proves the claim.
\end{proof}
Note that the fixed point equation on the level of the coefficients
$\vc{a}$ is always satisfied for $a_0$ and $a_1$.  If $k>1$, one
obtains the recursion
\begin{displaymath}
  a_k = \frac{1}{k-1} \, \sum_{n=1}^{k-1} a_n a_{k-n} \es,
\end{displaymath}
which shows that at most one fixed point with given mean can exist.

Let us now consider the case of discrete time first.  Analogously to
\eref{urdisctime}, define $\vc{a}(t) = \vc{a}(\vc{p}(t))$ as the
coefficients belonging to $\vc{p}(t)$, which are assumed to exist.  It
is clear from \eref{urgrecopsi}, \eref{ur1recopsi} and \eref{ur1recoa}
that $\vc{a}(t+1) = \Ra_1(\vc{a}(t))$.  We then have the following two
propositions.
\begin{proposition}
  \label{prop:ur1convdisc}
  Assume\/ $\vc{a}(0)$ to exist.  Then, in discrete time,
  \begin{displaymath}
    \lim_{t\to\infty} a_k(t) = \alpha^k
    \qquad\text{for all $k\ge0$.}
  \end{displaymath}
\end{proposition}
This result indicates that a weaker condition than the one of
\thref{ur1conv} may be sufficient for convergence of $\vc{p}(t)$.
\begin{proof}
  Clearly, $a_0(t) \equiv 1$, $a_1(t) \equiv \alpha$.  Furthermore, by
  the assumption and \eref{ur1recopsi}, the coefficients $a_k(t)$
  exist for all $k, t \in \NNzer$.  Now, assume that the claim holds
  for all $k \le n$ with some $n$ and let $k=n+1$.  According to the
  general properties of $\limsup$ and $\liminf$, we then have
  \begin{align*}
    \limsup_{t\to\infty} \, a_k(t+1) 
    &\le \frac{1}{k+1} \sum_{\ell=0}^k \limsup_{t\to\infty}
      \bigl( a_\ell(t) a_{k-\ell}(t) \bigr) \\
    &= \frac{k-1}{k+1} \alpha^k + \frac{2}{k+1} \limsup_{t\to\infty}\,a_k(t)
  \end{align*}
  and analogously with $\ge$ for $\liminf$.  This leads to
  \begin{displaymath}
    \frac{k-1}{k+1} \limsup_{t\to\infty}\,a_k(t) \le \frac{k-1}{k+1} \alpha^k
    \le \frac{k-1}{k+1} \liminf_{t\to\infty} a_k(t) \es,
  \end{displaymath}
  from which the claim follows for all $k \le n+1$ and, by induction
  over $n$, for all $k\ge0$.
\end{proof}
\begin{proposition}
  \label{prop:ur1contract}
  The recombinator\/ $\Ra_1$, acting on\/ $\X_{\alpha,\delta}$, is a
  strict contraction \wrt the metric\/ $d$ from \eref{urgmetric}, \ie
  there is a $C<1$ such that, for all elements $\vc{a}, \vc{b} \in
  \X_{\alpha,\delta}$,
  \begin{displaymath}
    d(\Ra_1(\vc{a}),\Ra_1(\vc{b})) \le C \, d(\vc{a},\vc{b}) \es.
  \end{displaymath}
\end{proposition}
\begin{proof}
  First consider, for $k\ge2$, without using the special choice of the
  $d_k$,
  \begin{equation}
    \label{eq:ur1contractineq}
    \begin{aligned}
      d(\Ra_1(\vc{a}),\Ra_1(\vc{b})) &=
      \sum_{k\ge2} d_k \frac{1}{k+1} \Bigl| 
        \sum_{\ell=0}^k (a_\ell\,a_{k-\ell} - b_\ell\,b_{k-\ell}) \Bigr| \\
      &= \sum_{k\ge2} d_k \frac{1}{k+1} \Bigl|
        \sum_{\ell=0}^k (a_\ell-b_\ell) (a_{k-\ell}+b_{k-\ell}) \Bigr| \\
      &\le \sum_{k\ge2} d_k \frac{2}{k+1} 
        \sum_{\ell=2}^k \delta^{k-\ell} |a_\ell-b_\ell| \\
      &= \sum_{\ell\ge2} d_\ell |a_\ell-b_\ell| 
        \sum_{k\ge\ell} \frac{2}{k+1} \delta^{k-\ell} \frac{d_k}{d_\ell} \es.
    \end{aligned}
  \end{equation}
  With the choice $d_k = (\gamma/\delta)^k$, where we had $\gamma <
  \frac13$, we can find, for $\ell\ge2$, an upper bound for the inner
  sum,
  \begin{displaymath}
    \sum_{k\ge\ell} \frac{2}{k+1} \delta^{k-\ell} \frac{d_k}{d_\ell} \le
    \frac23 \sum_{k\ge\ell} \gamma^{k-\ell} = \frac{2}{3-3\gamma} =: 
    C < 1 \es,
  \end{displaymath}
  which, together with \eref{ur1contractineq}, proves the claim.
\end{proof}
Together with Banach's fixed point theorem (compare \cite[\thm
V.18]{ReSi}), the two propositions imply that $\vc{a}(t)$ converges to
$(1,\alpha,\alpha^2,\ldots)$ \wrt the metric $d$, and that convergence
is exponentially fast.

In continuous time, we consider the time derivative of $\vc{a}(t) :=
\vc{a}(\vc{p}(t))$, which is, by \eref{urga},
\begin{equation}
  \label{eq:ur1ivpa}
  \tfrac{\dd{}}{\dd{t}} \vc{a}(t) =
  \tfrac{\dd{}}{\dd{t}} \vc{a}\bigl(\vc{p}(t)\bigr) =
  \vc{a}\bigl(\R_1(\vc{p}(t)) - \vc{p}(t)\bigr) =
  \Ra_1(\vc{a}(t)) - \vc{a}(t) \es.
\end{equation}
The following lemma ensures, together with \cite[\thm 7.6 and \remk
7.10(b)]{Ama}, that this initial value problem has a unique solution
for all $\vc{a}(0) = \vc{a}_0 \in \X_{\alpha,\delta}$.
\begin{lemma}
  \label{lem:ur1lipschitz}
  Consider the Banach space\/ $H_{\gamma/\delta}$ from \eref{hdelta},
  with some $0<\gamma<\tfrac13$, and its open subset\/ $Y = \{ \vc{x} \in
  H_{\gamma/\delta} : |x_k| < (2\delta)^k \}$.  Then, the
  recombinator\/ $\Ra_1$ from \eref{ur1recoa} maps\/ $Y$ into itself,
  satisfies a global Lipschitz condition, and is bounded on\/ $Y.$
  Furthermore, it is infinitely differentiable,\/ $\Ra_1 \in
  \C^\infty(Y,Y)$.
\end{lemma}
\begin{proof}
  For $\vc{x} \in Y,$ one has $|x_k| < (2\delta)^k$, hence
  $|\Ra_1(\vc{x})_k^{}| < (2\delta)^k$, with a similar argument as in
  the proof of \lref{ur1cont}.  Consequently, $\Ra_1(Y) \subset Y.$
  So, let $\vc{x}, \vc{y} \in Y.$  Then, similarly to the proof of
  \pref{ur1contract}, one shows the Lipschitz condition
  \begin{equation*}
    \|\Ra_1(\vc{x}) - \Ra_1(\vc{y})\| \le
    \sum_{\ell\ge0} \left(\frac\gamma\delta\right)^\ell |x_\ell - y_\ell|
      \sum_{k\ge\ell} \frac{2}{k+1} (2\gamma)^{k-\ell} \le
    \frac{2}{1-2\gamma} \|\vc{x}-\vc{y}\|
  \end{equation*}
  and, since $\|\vc{x}\| < 1/(1-2\gamma)$ in $Y$, the boundedness,
  \begin{displaymath}
    \|\Ra_1(\vc{x})\| \le \frac{1}{1-2\gamma} \|\vc{x}\| <
    \frac{1}{(1-2\gamma)^2} \es.
  \end{displaymath}
  \Wrt differentiability, consider, for sufficiently small $\vc{h} \in
  Y,$
  \begin{displaymath}
    \Ra_1(\vc{x}+\vc{h})_k^{} = 
    \Ra_1(\vc{x})_k^{} + \frac{2}{k+1} \sum_{\ell=0}^k x_{k-\ell}\,h_\ell +
      \Ra_1(\vc{h})_k^{} \es.
  \end{displaymath}
  Since
  \begin{align*}
    \|\Ra_1(\vc{h})\| &\le
    \sum_{k\ge0} \left(\frac\gamma\delta\right)^k \frac{1}{k+1}
      \sum_{\ell=0}^k |h_{k-\ell}| \, |h_\ell| \\
    &= \sum_{\ell\ge0} \left(\frac\gamma\delta\right)^\ell |h_\ell|
      \sum_{k\ge\ell} \left(\frac\gamma\delta\right)^{k-\ell}
        \frac{|h_{k-\ell}|}{k+1} \le \|\vc{h}\|^2 \es,
  \end{align*}
  it is clear that $\Ra_1$ is differentiable with linear (and thus
  continuous) derivative, whose Jacobi matrix is explicitly
  $\Ra_1'(\vc{x})_{k\ell}^{} = \frac{\partial}{\partial x_\ell}
  \Ra_1(\vc{x})_k^{} = \frac{2}{k+1} x_{k-\ell}$ if $k\ge\ell$ and
  zero otherwise, hence one has $\Ra_1 \in \C^1(Y,Y)$.  It is now
  trivial to show that $\Ra_1 \in \C^2(Y,Y)$ with constant second
  derivative and thus $\Ra_1 \in \C^\infty(Y,Y)$.
\end{proof}
\begin{proposition}
  \label{prop:ur1convcont}
  If\/ $\vc{a}_0 \in \X_{\alpha,\delta}$ for some\/ $\alpha, \delta$,
  then\/ $\vc{a}(t) \in \X_{\alpha,\delta}$ for all\/ $t\ge0$ and\/
  $\lim_{t\to\infty} d(\vc{a}(t), \vc{\alpha}) = 0$ with\/
  $\vc{\alpha} = (1,\alpha,\alpha^2,\alpha^3,\ldots)$.
\end{proposition}
\begin{proof}
  The first statement follows from \cite[\thm VI.2.1]{Mar} (see also
  \cite[\thm 16.5]{Ama}) since, due to the convexity of
  $\X_{\alpha,\delta}$, we have $\vc{a} + t (\Ra(\vc{a})-\vc{a}) \in
  \X_{\alpha,\delta}$ for every $\vc{a} \in \X_{\alpha,\delta}$ and
  $t\in\ccint{0,1}$, hence a subtangent condition is satisfied.  For
  the second, observe that $\Ra_1(\vc{\alpha}) = \vc{\alpha}$.  We now
  show that
  \begin{equation}
    \label{eq:ur1lyapunov}
    L(\vc{a}_0) = d(\vc{a}_0, \vc{\alpha})
  \end{equation}
  is a Lyapunov function, \cf \dref{lyapunov}.  With the notation of
  \lref{ur1lipschitz}, note that the compact metric space
  $\X_{\alpha,\delta}$ is contained in the open subset $Y$ of the
  Banach space $H_{\gamma/\delta}$.  The continuity of $L$ is obvious.
  Now, let $\vc{a}_0 \in \X_{\alpha,\delta}$ be given.  By
  \lref{ur1lipschitz} and \cite[\thm 9.5 and \remk 9.6(b)]{Ama}, the
  solution $\vc{a}(t)$ of \eref{ur1ivpa} is infinitely differentiable.
  Thus, for $t\in\ccint{0,1}$,
  \begin{equation}
    \label{eq:ur1lyapuineq}
    \begin{aligned}
      L(\vc{a}(t)) - L(\vc{a}_0) &=
      \|\vc{a}_0 + t (\Ra_1(\vc{a}_0)-\vc{a}_0) + \ord(t) - 
          \vc{\alpha}\| -
        \|\vc{a}_0 - \vc{\alpha}\| \\
      &\le t \bigl( \|\Ra_1(\vc{a}_0) - \Ra_1(\vc{\alpha})\| -
        \|\vc{a}_0 - \vc{\alpha}\| \bigr) + \ord(t) \es,
    \end{aligned}
  \end{equation}
  where $\ord(t)$ is the usual Landau symbol and represents some
  function that vanishes faster than $t$ as $t\to0$.  From this, by
  the strict contraction property of $\Ra_1$ (\pref{ur1contract}), the
  Lyapunov property \eref{dotl} follows, with equality if and only if
  $\vc{a}_0 = \vc{\alpha}$.  Since $\X_{\alpha,\delta}$ is compact,
  \thref{lyapunov} implies the claim.
\end{proof}

\pagebreak[2]
We are now able to give the previously postponed
\begin{theopargself}
  \begin{proof}[of \thref{ur1conv}]
    By \pref{urgainxad}, we have $\vc{a}(0) = \vc{a}(\vc{p}(0)) \in
    \X_{\alpha,\delta}$ with $\alpha = \frac12 m$ and some $\delta$.
    In discrete time, according to \propositions
    \ref{prop:ur1convdisc} and \ref{prop:ur1contract} and Banach's
    fixed point theorem (compare \cite[\thm V.18]{ReSi}), it then
    follows that $\vc{a}(t) \to \vc{\alpha} =
    (1,\alpha,\alpha^2,\ldots)$ \wrt the metric $d$.  Inserting
    \eref{ur1fixed} into \eref{urga} and letting $x = m/(m+2)$ yields
    \begin{align*}
      a_k^{} &=
      \sum_{\ell \ge k} \frac{\ell!}{(\ell-k)!(k+1)!} 
        (1-x)^2 (\ell+1) x^\ell \\
      &= (1-x)^2 \sum_{\ell \ge k} \binom{\ell+1}{k+1} x^\ell =
      \left(\frac{x}{1-x}\right)^k = \alpha^k .
    \end{align*}
    The claim now follows from \lref{urgpnormconv}.  Similarly, in
    continuous time, the claim follows from \pref{ur1convcont}.
  \end{proof}
\end{theopargself}

Let us finally note
\begin{proposition}
  For the dynamics described by \eref{ur1ivpa}, the fixed point\/
  $\vc{\alpha}$ from \pref{ur1convcont} is exponentially stable.
\end{proposition}
\begin{proof}
  Let $\vc{a}_0 \in \X_{\alpha,\delta}$ be arbitrary.  The Lyapunov
  function from the proof of \pref{ur1convcont} satisfies, as a
  consequence of \eref{ur1lyapuineq} and \pref{ur1contract},
  \begin{displaymath}
    \dot{L}(\vc{a}_0) \le 
    d(\Ra_1(\vc{a}_0), \Ra_1(\vc{\alpha})) - d(\vc{a}_0, \vc{\alpha}) \le
    - (1-C) \, d(\vc{a}_0, \vc{\alpha}) \es,
  \end{displaymath}
  with $0 < C < 1$.  From this, together with \eref{ur1lyapunov} and
  \cite[\thm 18.7]{Ama}, the claim follows.
\end{proof}

\noindent
\textbf{Remark.}
In a related UC model introduced by Takahata \cite{Tak}, for which
\begin{displaymath}
  T_{ij,k\ell} = \delta_{i+j,k+\ell} \frac{1}{k+\ell+1} \es,
\end{displaymath}
the recombinator $\Ra_1$ appears for the coefficients
$\vc{b}(\vc{p})_k^{} = (k+1) \, \vc{a}(\vc{p})_k^{}$, where
$\vc{b}(\vc{p})_1^{}$ is the mean copy number $m$.  The above results
then imply, under the appropriate condition on $\vc{p}(0)$, that
$\vc{b}(t) \to (1,m,m^2,\ldots)$ as $t\to\infty$ both in discrete and
in continuous time.  This corresponds to convergence of $\vc{p}(t)$ to
the fixed point $\vc{p}$ with $p_k^{} = \frac{1}{m+1}
(\frac{m}{m+1})^k$.

\section{The intermediate parameter regime}
\label{sec:urqgen}

In this section, $q$ may take any value in $\ccint{0,1}$.  \Wrt
reversibility of fixed points, one finds
\begin{proposition}
  For parameter values\/ $q\in\ooint{0,1}$, any fixed point\/ $\vc{p}
  \in \M_1^+$ of the recombinator\/ $\R_q$, given by \eref{reco} and
  \eref{urtq}, satisfies\/ $p_k^{}>0$ for all\/ $k\ge0$ \pr{unless it
    is the trivial fixed point\/ $\vc{p} = (1,0,0,\ldots)$ we
    excluded}.  None of these extra fixed points is reversible.
\end{proposition}
\begin{proof}
  Let a non-trivial fixed point $\vc{p}$ be given and choose any $n>0$
  with $p_n^{}>0$.  Observe that $T^{(q)}_{n+1 \is n-1, nn} > 0$ for
  $0 < q < 1$ and hence
  \begin{displaymath}
    p_{n\pm1}^{\vp} = \R_q(\vc{p})_{n\pm1}^{\vp} 
    = \sum_{j,k,\ell\ge0} T^{(q)}_{n\pm1 \is j,k\ell} \, 
      p_k^{\vp} \, p_\ell^{\vp}
    \ge T^{(q)}_{n+1 \is n-1, nn} \, p_n^{\vp} \, p_n^{\vp} > 0 \es.
  \end{displaymath}
  The first statement follows now by induction.  

  For the second statement,
  evaluate the reversibility condition \eref{urrev} for all
  combinations of $i$, $j$, $k$, $\ell$ with $i+j = k+\ell \le 4$.
  This leads to four independent equations.  Three of them can be
  transformed to the recursion
  \begin{displaymath}
    p_k^{} 
    = \frac{(k+1)q}{2(k-1)+2q} \frac{p_1^{}}{p_0^{}} \, p_{k-1}^{} \es,
    \qquad k \in \{2,3,4\} \es,
  \end{displaymath}
  from which one derives explicit equations for all $p_k^{}$ with $k
  \in \{2,3,4\}$ in terms of $p_0^{}$ and $p_1^{}$.  Inserting the one
  for $p_2^{}$ into the remaining equation yields another equation for
  $p_4^{}$ in terms of $p_0^{}$ and $p_1^{}$, which contradicts the
  first equation for all $q \in \ooint{0,1}$, as is easily verified.
\end{proof}
So, non-trivial fixed points for $0 < q < 1$ are not reversible, and
thus much more difficult to determine.  Our most general result so far
is
\begin{theorem}
  \label{thm:urqfixed}
  If\/ $\vc{p}(0) \in P_{\alpha,\delta}$ for some\/ $\alpha, \delta$,
  then\/ $\vc{p}(t) \in P_{\alpha,\delta}$ for all times\/
  $t\in\NNzer$, respectively\/ $t\in\RRnn$, and\/ $\R_q$ has a fixed
  point in\/ $P_{\alpha,\delta}$.
\end{theorem}
The proof is based on the fact that $\R_q$ is, in a certain sense,
monotonic in the parameter $q$.  This is stated in
\begin{proposition}
  \label{prop:urqainxad}
  Assume\/ $\vc{a}(\vc{p}) \in \X_{\alpha,\delta}$ for some\/
  $\alpha$, $\delta$.  Then, \wrt the partial order introduced before
  \pref{urgpadcomp}, $\vc{a}(\R_q(\vc{p})) \le
  \vc{a}(\R_{q'}(\vc{p}))$ for all\/ $0 \le q \le q' \le 1$.  In
  particular,\/ $\vc{a}(\R_q(\vc{p})) \in \X_{\alpha,\delta}$ for
  all\/ $0 \le q \le 1$.
\end{proposition}
To show this, we need three rather technical lemmas.  The first one
collects formal conditions on the difference of two discrete probability 
distributions
$T^{(q)}_{ij,k\ell}$ with different parameter values (but $j=k+\ell-i$
and the same fixed $k$, $\ell$).  These are then verified in our case.
\begin{lemma}
  \label{lem:urqineq}
  Let the numbers\/ $x_i\in\RR$ \pr{$0 \le i \le r$ with some\/
    $r\in\NNzer$} satisfy the following three conditions:
  \begin{gather}
    \label{eq:urqineqcond1}
    \sum_{i=0}^r x_i = 0 \es.\\
    \label{eq:urqineqcond2}
    x_{r-i} = x_i \quad\text{for all\/ $0 \le i \le r$.}\\
    \label{eq:urqineqcond3}
    \text{There is an integer\/ $n$ such that }
    \begin{cases} x_i \ge0 \; : \; 0 \le i \le n \\
      x_i<0 \; : \; n < i \le \lfloor \frac{r}{2} \rfloor \end{cases}
    \!\!\!\!\!.
  \end{gather}
  Further, let\/ $f_i\in\RR$ \pr{$0 \le i \le r$} be given with
  \begin{equation}
    \label{eq:urqineqcond4}
    0 \le f_1 - f_0 \le f_2 - f_1 \le \ldots \le f_r - f_{r-1} \es.
  \end{equation}
  Then, we have
  \begin{displaymath}
    \sum_{i=0}^r f_i x_i \ge 0 \es.
  \end{displaymath}
\end{lemma}
\begin{proof}
  Let us first consider the trivial cases.  If $x_i\equiv0$,
  everything is clear, so let $x_i\not\equiv0$.  If $r\le1$ then
  $x_i\equiv0$, so let $r\ge2$, and thus $n\le\frac{r}{2}-1$.  Define
  $x_{\frac{r}{2}} = f_{\frac{r}{2}} = 0$ for odd $r$.  Then, we can
  write
  \begin{equation}
    \label{eq:urqineq1}
    \sum_{i=0}^r f_i x_i = \sum_{i=0}^n \left( f_i + f_{r-i} \right) x_i +
    \sum_{i=n+1}^{\lceil\frac{r}{2}\rceil-1} \left( f_i + f_{r-i}
      \right) x_i + f_{\frac{r}{2}} x_{\frac{r}{2}} \es.
  \end{equation}
  Furthermore, for $r-i \ge i$, due to \eref{urqineqcond4},
  \begin{displaymath}
    f_i + f_{r-i} = f_{i-1} + f_{r-i+1} + (f_i - f_{i-1}) - 
    (f_{r-i+1} - f_{r-i}) \le f_{i-1} + f_{r-i+1} \es.
  \end{displaymath}
  Now, define $C := \sum_{i=0}^n x_i = -
  \sum_{i=n+1}^{\lceil\frac{r}{2}\rceil-1} x_i - \tfrac12
  x_{\frac{r}{2}} > 0$, and the claim follows with \eref{urqineq1},
  since $r-n \ge n+1$ by assumption:
  \begin{multline*}
    \sum_{i=0}^r f_i x_i \ge 
    C \left[ f_n + f_{r-n} - f_{n+1} - f_{r-n-1} \right] \\
    = C \left[ (f_{r-n} - f_{r-n-1}) - (f_{n+1} - f_n) \right] \ge 0 \es.
    \qedhere
  \end{multline*}
\end{proof}
\begin{lemma}
  \label{lem:urqcoeff}
  Let\/ $j\in\NNzer$ be fixed and\/ $f_i = (i)_j^{}$,
  $i\in\NNzer$, where\/ $(i)_j$ is the falling factorial, which
  equals\/ $1$ for\/ $j=0$ and\/ $i(i-1)\cdots(i-j+1)$ for\/
  $j>0$, hence\/ $\frac{i!}{(i-j)!}$ for\/ $i \ge j$.  Then
  condition \eref{urqineqcond4} is satisfied.
\end{lemma}
\begin{proof}
  For $j=0$, condition \eref{urqineqcond4} is trivially true.
  Otherwise, each $f_i$ is a polynomial of degree $j$ in $i$ with
  zeros $\{0,1,\ldots,j-1\}$, hence we have the equality $0 = f_1 -
  f_0 = \ldots = f_{j-1} - f_{j-2}$.  Then, for $i \ge j-1$, the
  polynomial and all its derivatives are increasing functions since
  $\lim_{i\to\infty} f_i = \infty$.  Therefore, for $i \ge j-1$, we
  have $0 \le f_{i+1} - f_i \le f_{i+2} - f_{i+1}$.  Hence
  \eref{urqineqcond4} holds.
\end{proof}
\begin{lemma}
  \label{lem:urqineqcond123}
  For\/ $0 \le q \le q' \le 1$ and all\/ $k, \ell$, Equations
  \eref{urqineqcond1}--\eref{urqineqcond3} are true for\/ $r=k+\ell$
  and $x_i = T^{(q')}_{ik\ell} - T^{(q)}_{ik\ell}$, where\/
  $T^{(q)}_{ik\ell} = T^{(q)}_{ij,k\ell}$ with\/ $j = k+\ell-i$.
\end{lemma}
\begin{proof}
  The validity of \eref{urqineqcond1} and \eref{urqineqcond2} is clear
  from the normalization \eref{sumt} and the symmetry of the
  $T^{(q)}_{ik\ell}$.  For \eref{urqineqcond3}, let $k \le \ell$
  without loss of generality.  In the trivial cases $q=q'$ or $k=0$,
  choose $n=\lfloor\frac{r}{2}\rfloor$.  Otherwise, $x_i =
  T^{(q')}_{ik\ell} - T^{(q)}_{ik\ell} < 0$ for $k \le i \le
  \lfloor\frac{r}{2}\rfloor$, since $C^{(q')}_{k\ell} <
  C^{(q)}_{k\ell}$, and $x_0 > 0$.  For $0 \le i \le k$, consider
  \begin{displaymath}
    y_i = \frac{x_i}{T^{(q)}_{ik\ell}} + 1
    = \frac{C^{(q')}_{k\ell}}{C^{(q)}_{k\ell}} 
      \left(\frac{q'}{q}\right)^{k-i} \es.
  \end{displaymath}
  Here, the first factor is less than 1, the second is equal to 1 for
  $k=i$, greater than 1 for $0 \le k < i$, and strictly decreasing
  with $i$.  Since $x_i\ge0$ if and only if $y_i\ge1$, there is an
  index $n$ with the properties needed.
\end{proof}
\begin{theopargself}
  \begin{proof}[of \pref{urqainxad}]
    We assume $0\le q\le q' \le 1$.
    Lemmas \ref{lem:urqineq}--\ref{lem:urqineqcond123} imply, for all
    $k, \ell, j\in\NNzer$ with $k+\ell \ge j$,
    \begin{displaymath}
      \sum_{i=j}^{k+\ell} \tfrac{i!}{(i-j)!} T^{(q)}_{ik\ell} \le 
      \sum_{i=j}^{k+\ell} \tfrac{i!}{(i-j)!} T^{(q')}_{ik\ell} \es.
    \end{displaymath}
    Then, since $T^{(q)}_{ik\ell} = 0$ for $i >k +\ell$,
    \begin{align*}
      \vc{a}(\R_q(\vc{p}))_j^{} &= \tfrac{1}{(j+1)!} \sum_{i \ge j}
      \tfrac{i!}{(i-j)!} \R_q(\vc{p})_i^{} = 
      \tfrac{1}{(j+1)!} \sum_{i \ge j} \tfrac{i!}{(i-j)!}
      \sum_{k,\ell\ge0} T^{(q)}_{ik\ell} p_k^{\vp} p_\ell^{\vp} \\
      &=\tfrac{1}{(j+1)!} \mskip-4mu\sum_{k,\ell\ge0}\mskip-4mu 
      p_k^{\vp} p_\ell^{\vp} \sum_{i \ge j}
      \tfrac{i!}{(i-j)!} T^{(q)}_{ik\ell} \\
      &\le \tfrac{1}{(j+1)!} \mskip-4mu\sum_{k,\ell\ge0}\mskip-4mu 
      p_k^{\vp} p_\ell^{\vp} \sum_{i \ge j}
      \tfrac{i!}{(i-j)!} T^{(q')}_{ik\ell}
      = \vc{a}(\R_{q'}(\vc{p}))_j^{} \es.
    \end{align*}
    From this, together with \lref{ur1cont}, the claim follows.
  \end{proof}
\end{theopargself}
\begin{theopargself}
  \begin{proof}[of \thref{urqfixed}]
    According to \pref{urqainxad}, $\R_q$ maps $P_{\alpha,\delta}$
    into itself, and thus, in discrete time, $\vc{p}(t) \in
    P_{\alpha,\delta}$ for every $t\in\NNzer$.  The analogous
    statement is true for continuous time $t\in\RRnn$.  To see this,
    consider $P_{\alpha,\delta}$ as a closed subset of $\Ll{1}$.
    Recall that $\R_q - \Id$ is globally Lipschitz on $\Ll{1}$ by
    \pref{recolipschitz}.  Moreover, for any $\vc{p} \in
    P_{\alpha,\delta}$ and $t\in\ccint{0,1}$, \pref{urgpadcomp} tells
    us that
    \begin{displaymath}
      \vc{p} + t (\R_q(\vc{p}) - \vc{p}) = 
      (1-t) \vc{p} + t \R_q(\vc{p}) \in P_{\alpha,\delta} \es.
    \end{displaymath}
    This implies the positive invariance of $P_{\alpha,\delta}$ by
    \cite[\thm VI.2.1]{Mar} (see also \cite[\thm 16.5]{Ama}).  The
    existence of a fixed point once again follows from the
    Leray--Schauder--Tychonov theorem \cite[\thm V.19]{ReSi}.
  \end{proof}
\end{theopargself}

On the basis of the above analysis, and further numerical work done
to investigate the fixed point properties \cite{Red,ShAt},
it is plausible that, given the mean copy number $m$, never more than
one fixed point for $\R_q$ exists.  Due to the global convergence
results at $q=0$ and $q=1$, any non-uniqueness in the vicinity of
these parameter values could only come from a bifurcation, not from an
independent source.  Numerical investigations indicate that no
bifurcation is present, but this needs to be analyzed further.

Furthermore, the Lipschitz constant for $\Ra_q$ can be expected to be
continuous in the parameter $q$, hence to remain strictly less than 1
on the sets $X_{\alpha,\delta}$ in a neighborhood of $q=1$.  So, at
least locally, the contraction property should be preserved.
Nevertheless, we do not expand on this here since it seems possible to
use a rather different approach \cite{Hof03}, which has been used for
similar problems in game theory, to establish a slightly weaker type
of convergence result for all $0 < q < 1$, and probably even on the
larger compact set $\M_{1,m,C}^+$ of \lref{ur0dp}.

\section{Concluding remarks}
\label{sec:remarks}

In this article, we have shown that, for the extreme parameter values $q=0$
(internal UC) and $q=1$ (random UC), any initial configuration satisfying a
specific condition converges to one of the known fixed points, both in
discrete and continuous time.  The condition to be met is, for $q=0$, the
existence of the $r$-th moment ($r>1$, see \thref{ur0conv}), respectively, for
$q=1$, that the corresponding generating function has a radius of
convergence $\rho>1$ (\thref{ur1conv}).  
Convergence takes place in the total variation
norm in all cases.  As argued in the previous section, similar results can be
expected for the intermediate parameter values as well.

% We have seen in the preceding sections that, definitely for the extreme cases
% $q=0$ and $q=1$ and presumably for the intermediate values as well, the
% deterministic dynamics converges to a unique equilibrium solution in both
% discrete and continuous time.  This correspond
These results are valid for deterministic dynamics and thus correspond
to the case of infinite populations.  \Wrt biological
relevance, however, we add some arguments that it is reasonable to expect this
to be a good description for large but finite populations as well, \ie for the
underlying (multitype) branching process.  For finite state spaces, such as in
the mutation--selection models discussed in \cite{HRWB}, the results by Ethier
and Kurtz \cite[\thm 11.2.1]{EtKu86} and the generalization \cite[\thm
V.7.2]{AtNe72} of the Kesten--Stigum theorem \cite{KeSt66,KLPP97} guarantee
that in the infinite population limit the relative genotype frequencies of the
branching process converge almost surely to the deterministic solution (if the
population does not go to extinction).  Since for the UC models considered
here the equilibrium distributions are exponentially small for large copy
numbers (owing to \thref{urqfixed} also for $q\in\ooint{0,1}$), one can expect
these systems to behave very much like ones with finitely many genotypes.
This is also supported by several simulations.  Nevertheless, this questions
deserves further attention.

\begin{acknowledgement}
  We thank Ellen Baake, Joachim Hermisson, and Josef Hofbauer for
  their cooperation, and Odo Diekmann and Anton Wakolbinger for
  clarifying discussions.  O.R. gratefully acknowledges support by a
  PhD scholarship of the Studienstiftung des deutschen Volkes.  The
  authors express their gratitude to the Erwin Schr\"odinger
  International Institute for Mathematical Physics in Vienna for
  support during a stay in December 2002, where the manuscript for
  this article was completed.
\end{acknowledgement}

\end{document}